\documentclass[reqno]{myart}

\usepackage{mathtools}

\usepackage{stmaryrd}
\usepackage[
  warnings-off={mathtools-colon,mathtools-overbracket}]{unicode-math}
\let\dfn\emph
\setoperatorfont{\textup}

\AfterEndPreamble{%
  \let\phi\varphi%
  \let\epsilon\varepsilon%
  \let\geq\geqslant%
  \let\leq\leqslant%
  \let\div\relax%
  \newcommand{\div}{\mathord\divslash}%
  \renewcommand{\div}{/}%
  \renewcommand{\setminus}{\char"2216}
  \renewcommand{\emptyset}{\char"2300}
  \let\sharp\mathoctothorpe
}
\newcommand{\defeq}{\coloneq}
\newcommand{\blank}{\mathord{-}}
\NewDocumentCommand\pair{O{\blank, \blank}}{\langle#1\rangle}
\newcommand{\abs}[1]{\vert#1\vert}

\newcommand{\tb}{\textup{b}}
\newcommand{\te}{\textup{e}}
\newcommand{\uD}{\symup{D}}

\newcommand{\uR}{\symup{R}}
\newcommand{\uM}{\symup{M}}
\newcommand{\Ab}{\symup{Ab}}
\newcommand{\kk}{\symbb{k}}
\newcommand{\NN}{\symbb{N}}
\newcommand{\ZZ}{\symbb{Z}}

\newcommand{\cC}{\symcal{C}}

\newcommand{\cM}{\symcal{M}}
\newcommand{\cN}{\symcal{N}}
\newcommand{\cP}{\symcal{P}}

\newcommand{\cT}{\symcal{T}}

\newcommand{\id}{\textup{id}}
\newcommand{\op}{\textup{op}}
\newcommand{\opp}[1]{{#1}^{\op}}
\newcommand{\Ao}{\opp{A}}
\newcommand{\Ae}{A^{\te}}

\DeclareMathOperator{\coker}{coker}
\DeclareMathOperator{\Kdim}{Kdim}
\DeclareMathOperator{\im}{im}
\let\lim\relax
\DeclareMathOperator*{\lim}{lim}

\DeclareMathOperator{\Hom}{Hom}
\DeclareMathOperator{\gHom}{\underline{Hom}}
\DeclareMathOperator{\End}{End}
\DeclareMathOperator{\gEnd}{\underline{End}}
\DeclareMathOperator{\AutOp}{Aut}
\newcommand{\Aut}[1]{\AutOp(#1)}
\DeclareMathOperator{\GrAutOp}{GrAut}
\newcommand{\GrAut}[1]{\GrAutOp(#1)}
\DeclareMathOperator{\Ext}{Ext}
\DeclareMathOperator{\gExt}{\underline{Ext}}
\DeclareMathOperator{\Tor}{Tor}
\DeclareMathOperator{\RHom}{RHom}
\DeclareMathOperator{\RgHom}{\uR\underline{Hom}}
\DeclareMathOperator{\GKdim}{GKdim}
\DeclareMathOperator{\add}{add}

\DeclareMathOperator{\ind}{ind}
\DeclareMathOperator{\PGL}{PGL}
\DeclareMathOperator{\charOp}{char}
\DeclareMathOperator{\simOp}{sim}
\newcommand{\simp}[1]{\simOp(#1)}

\DeclareMathOperator{\ProjOp}{Proj}
\newcommand{\Proj}[1]{\ProjOp(#1)}
\DeclareMathOperator{\QGrOp}{QGr}
\newcommand{\QGr}[1]{\QGrOp(#1)}
\DeclareMathOperator{\qgrOp}{qgr}
\newcommand{\qgr}[1]{\qgrOp(#1)}
\DeclareMathOperator{\cohOp}{coh}
\newcommand{\coh}[1]{\cohOp(#1)}

\DeclareMathOperator{\rankOp}{rank}
\newcommand{\rank}[1]{\rankOp(#1)}
\DeclareMathOperator{\depthOp}{depth}
\let\dep\depthOp
\newcommand{\depth}[1]{\depthOp(#1)}
\DeclareMathOperator{\socOp}{soc}
\newcommand{\soc}[1]{\socOp(#1)}

\DeclareMathOperator{\idim}{inj. dim}
\DeclareMathOperator{\gldim}{gl. dim}

\DeclareMathOperator{\modOp}{mod}
\let\mod\relax
\newcommand{\mod}[1]{\modOp(#1)}
\DeclareMathOperator{\ModOp}{Mod}
\newcommand{\Mod}[1]{\ModOp(#1)}
\DeclareMathOperator{\GrOp}{Gr}
\newcommand{\Gr}[1]{\GrOp(#1)}
\DeclareMathOperator{\grOp}{gr}
\newcommand{\gr}[1]{\grOp(#1)}
\DeclareMathOperator{\TorsOp}{Tors}
\newcommand{\Tors}[1]{\TorsOp(#1)}
\DeclareMathOperator{\torsOp}{tors}
\newcommand{\tors}[1]{\torsOp(#1)}
\DeclareMathOperator{\CMOp}{CM}
\NewDocumentCommand{\CM}{sm}{%
  \IfBooleanTF#1%
    {\CMOp^0(#2)}
    {\CMOp(#2)}}
\DeclareMathOperator{\sCMOp}{\underline{CM}}
\NewDocumentCommand{\sCM}{sm}{%
  \IfBooleanTF#1%
    {\sCMOp^0(#2)}
    {\sCMOp(#2)}}
\DeclareMathOperator{\GL}{GL}

\newcommand{\Dbqgr}[1]{\uD^{\tb}(\qgr{#1})}

\newcommand{\Dbmod}[1]{\uD^{\tb}(\mod{#1})}
\newcommand{\Dbcoh}[1]{\uD^{\tb}(\coh{#1})}

\makeatletter
\newcommand{\raisemath}[1]{\mathpalette{\raisem@th{#1}}}
\newcommand{\raisem@th}[3]{\raisebox{#1}{$#2#3$}}
\makeatother

\usepackage[
  backend=biber,
  style=alphabetic,sorting=nyt,doi=false,url=false,
  giveninits=true,
  minbibnames=3,minalphanames=3,maxalphanames=4,maxbibnames=4,
  citetracker=true
  ]{biblatex}
\defbibheading{bibliography}[References]{%
  \section*{#1}}

\newcommand\myprintbibliography{%
  \singlespacing%
  \sloppy\hbadness 10000\relax%
    \printbibliography
  \fussy}

\usepackage{tikz-cd}
\usetikzlibrary{nfold}
\tikzcdset{
  arrow style=Latin Modern,
  diagrams={>={Straight Barb[scale=0.8]}}
}
\tikzset{
  no line/.style={draw=none,
    commutative diagrams/every label/.append style={/tikz/auto=false}},
  from/.style args={#1 to #2}{to path={(#1)--(#2)\tikztonodes}}}

\usepackage{enumitem}
\setlist{nosep}
\setlist{leftmargin=0pt,labelwidth=0pt,
  labelsep=\fontdimen2\font,itemindent=!,align=left}
\setlist[1]{labelindent=\parindent}
\setlist[enumerate,1]{label=(\arabic*)}
\setlist[2]{labelindent=2\parindent}
\setlist[enumerate,2]{label=(\roman*)}

\usepackage{varioref}
\usepackage[draft=false]{hyperref}
\usepackage[capitalise,nameinlink,noabbrev]{cleveref}
\usepackage{crossreftools}
\crefname{enumi}{}{}
\crefname{diagram}{diagram}{diagrams}
\crefname{equation}{}{}
\AtBeginDocument{\let\ref\relax\let\ref\Cref}
\makeatletter
\newcommand{\enumilabel}[2]{%
  \phantomsection
  #1\protected@edef\cref@currentlabel{[enumi][][]#1}\label[enumi]{#2}%
}
\makeatother

\usepackage{booktabs}
\usepackage{makecell}
\usepackage{multirow}
\usepackage{caption}
\makeatletter
\numberwithin{table}{section}
\AtEndPreamble{\renewcommand\c@table\c@theorem}
\makeatother

\usepackage{quiver}

\newtheorem{theorem}{Theorem}[section]
\newtheorem{lemma}[theorem]{Lemma}
\AddToHook{env/lemma/begin}{\crefalias{theorem}{lemma}}
\newtheorem{proposition}[theorem]{Proposition}
\AddToHook{env/proposition/begin}{\crefalias{theorem}{proposition}}
\newtheorem{corollary}[theorem]{Corollary}
\AddToHook{env/corollary/begin}{\crefalias{theorem}{corollary}}
\theoremstyle{definition}
\newtheorem{definition}[theorem]{Definition}
\AddToHook{env/definition/begin}{\crefalias{theorem}{definition}}
\newtheorem{example}[theorem]{Example}
\AddToHook{env/example/begin}{\crefalias{theorem}{example}}

\theoremstyle{remark}
\newtheorem{remark}[theorem]{Remark}
\AddToHook{env/remark/begin}{\crefalias{theorem}{remark}}
\numberwithin{equation}{section}

\title[Tilting and Cluster Tilting for Noncommutative Quadrics]
  {Tilting and Cluster Tilting for
  Noncommutative Isolated Quadric Hypersurface Singularities}

\author[Chen X.-Y.]{Xuyang Chen}
\address{School of Mathematical Sciences,
  Fudan University, Shanghai 200433, China}
\email{xuyangchen21@m.fudan.edu.cn}

\author[Li H.-N.]{Haonan Li}
\address{School of Mathematics,
  Shanghai University of Finance and Economics,
  Shanghai 200433, China}
\email{lihaonan@mail.shufe.edu.cn}
\thanks{Haonan Li is supported by the Fundamental Research Funds for the Central Universities and the National Natural Science Foundation of China (Grant No.\@ 1260011206).}

\author[Wu Q.-S.]{Quanshui Wu}
\address{School of Mathematical Sciences,
  Fudan University, Shanghai 200433, China}
\email{qswu@fudan.edu.cn}
\thanks{Quanshui Wu is supported by the National Natural Science Foundation of China (Grant No.\@ 12471032).}


\date{\today}


\begin{document}
  \begin{abstract}


We study tilting objects and cluster tilting modules for noncommutative isolated quadric hypersurface singularities $A$
over an algebraically closed field $\kk$ of characteristic different from $2$.
We prove that, for each $d\geq 2$, there exists a noncommutative isolated quadric hypersurface singularity of dimension $d$ admitting a graded $(d-1)$-cluster tilting module if and only if $d=2^n+1$ for some integer $n\geq 0$. 
Examples are constructed in every such dimension.
In contrast, no commutative isolated quadric hypersurface singularity of dimension $d\geq 4$ admits a graded $(d-1)$-cluster tilting module.
When the singularity $A$ has dimension at least $2$,
we also construct a tilting object in $\Dbqgr{A}$
by generalizing Kapranov's full strong exceptional sequence,
thereby obtaining an equivalence between $\Dbqgr{A}$ and the bounded derived category of a finite-dimensional algebra.
When $A$ has dimension $2$,
we use this tilting object to construct a noncommutative resolution of $A$
that is realized as a graded twist of a graded preprojective algebra of Euclidean type $\widetilde{A_1}$ or $\widetilde{D_4}$.
We establish a correspondence between the graded isomorphism classes of these twisted preprojective algebras
and the graded isomorphism classes of noncommutative isolated quadric surface singularities.
Consequently, we classify these singularities in terms of conjugacy classes in the groups $\PGL_2(\kk)$ and $S_4$.

  \end{abstract}

  \maketitle

  \setcounter{tocdepth}{1}
  \vspace{-2ex}
  \begingroup\small\tableofcontents\endgroup


\section{Introduction}\label{sec:intro}

Tilting and cluster tilting theory play important roles in
representation theory, commutative algebra, and algebraic geometry;
see Iyama's survey \cite{iyamaTiltingCohenMacaulayRepresentations2018}.
For example, Iyama and Takahashi \cite{iyaTiltingClusterTilting2013} used cluster tilting theory to
construct tilting objects in stable categories of
graded Cohen--Macaulay modules over commutative Gorenstein isolated quotient singularities, thereby showing that these stable categories are triangulated equivalent to derived categories of certain finite-dimensional algebras.
Mori and Ueyama \cite{morStableCategoriesGraded2016} extended this result to
noncommutative Gorenstein isolated quotient singularities arising from
finite ample group actions on noetherian Koszul Artin--Schelter regular algebras.
They constructed tilting objects in the stable categories of graded Cohen--Macaulay modules using tilting objects in $\qgr{A}$,
where $\qgr{A}$ denotes the quotient of the category of
finitely generated graded right $A$-modules modulo its full subcategory 
of finite-dimensional modules.
Following the fundamental work of Artin and Zhang \cite{artinNoncommutativeProjectiveSchemes1994}, $\qgr{A}$ is regarded as the category of coherent sheaves on the noncommutative projective scheme associated with $A$.


Affine cones over smooth projective quadrics
form another class of isolated singularities.
The spinor bundles on the quadrics arise from indecomposable
non-projective graded Cohen--Macaulay modules
over their homogeneous coordinate rings.
Kapranov \cite[Sec.\@ 4]{kapranovDerivedCategoriesCoherent1988} showed that these bundles, together with suitable twists
of the structure sheaf, form full strong exceptional sequences
in the bounded derived categories of coherent sheaves on the quadrics.

As analogues of projective quadrics,
noncommutative projective quadric hypersurfaces have been extensively studied in noncommutative projective geometry in recent years;
see, for example, \cite{smithNoncommutativeQuadricSurfaces2013,
moriNoncommutativeKnorrersPeriodicity2022,
heCliffordDeformationsKoszul2024,
huCliffordQuadraticComplete2024,
baumannCentralCurvesNoncommutative2026,
iyamaCohenMacaulayRepresentationsArtinSchelter2025,
liuSkewKnorrersPeriodicity2026,
moriIrreducibleNoncommutativeQuadrics2026,
}.
Smith and Van den Bergh \cite{smithNoncommutativeQuadricSurfaces2013}
established fundamental results on
noncommutative projective quadric surfaces and their
homogeneous coordinate algebras, which are quotients
of quantum polynomial algebras by central regular quadratic elements.
Mori and Ueyama \cite{moriNoncommutativeKnorrersPeriodicity2022}
extended several of these results
to quotients by normal regular quadratic elements, among other contributions.

In this paper, we work over an algebraically closed field $\kk$
with $\charOp(\kk)\neq 2$.
We study tilting objects, cluster tilting modules,
and noncommutative resolutions
for noncommutative isolated quadric hypersurface singularities.
Let $A=S/(f)$, where $S$ is a quantum polynomial algebra of dimension $d+1$
and $f\in S_2$ is normal and regular.
In the language of noncommutative projective geometry,
$A$ is a homogeneous coordinate algebra
of a noncommutative projective quadric hypersurface of dimension $d-1$.
We call $A$ a noncommutative quadric hypersurface singularity of dimension $d$.
Here, ``singularity" refers to the fact that $A$ has infinite global dimension.
We say that $A$ is a graded isolated singularity
if $\qgr{A}$ has finite global dimension.
In the commutative case, this condition means that the associated
projective scheme is smooth, while its affine cone has an isolated singularity at the vertex.
A noncommutative quadric hypersurface singularity $A$ with a graded isolated singularity is called
a noncommutative isolated quadric hypersurface singularity.
When $d=2$, $A$ is called a noncommutative isolated quadric surface singularity.


The notion of a noncommutative crepant resolution was introduced in
\cite{vandenberghNoncommutativeCrepantResolutions2004},
and was later shown to be closely connected with cluster tilting theory
\cite{iyamaAuslanderCorrespondence2007,
iyamaFominZelevinskyMutationTilting2008,iyamaMaximalModificationsAuslander2014}.
For noncommutative algebras, some related notions of resolution were introduced and studied in
\cite{qinNoncommutativeQuasiresolutions2019,hePreresolutionsNoncommutativeIsolated2022}.
For Artin--Schelter Gorenstein isolated singularities
of dimension $d\geq 2$,
Li, Shen, and Wu introduced the notion of a noncommutative resolution
(\cite[Def.\@ 6.1]{liNoncommutativeResolutionsASGorenstein2026})
and showed that noncommutative resolutions are
the graded endomorphism algebras of graded $(d-1)$-cluster tilting modules
\cite[Thm.\@ 6.6]{liNoncommutativeResolutionsASGorenstein2026}
(see \ref{dfn:NCCR,thm:NCCR-ct}).

In \ref{sec:tilting}, we consider tilting objects in $\Dbqgr{A}$
arising from full strong exceptional sequences. 
Iyama, Kimura, and Ueyama constructed
a full strong exceptional sequence in $\Dbqgr{A}$
for a noncommutative isolated quadric hypersurface singularity $A$
of dimension $d\geq 2$
\cite[Thm.\@ 9.8(3)]{iyamaCohenMacaulayRepresentationsArtinSchelter2025}.
To study noncommutative resolutions,
we construct an alternative full strong exceptional sequence
from degree shifts of $A$
and indecomposable non-projective graded Cohen--Macaulay $A$-modules
generated in degree 0.
Our sequence generalizes Kapranov's sequence for
commutative projective quadrics
\cite[Sec.\@ 4]{kapranovDerivedCategoriesCoherent1988}
and Ueyama's sequence for $(\pm 1)$-skew quadric hypersurface singularities
\cite[Thm.\@ 3.21]{ueyamaDerivedCategoriesSkew2023}.
Let $\pi\colon\gr{A}\to\qgr{A}$ denote the quotient functor.

\begin{theorem}[{\ref{thm:tilting-hypersurface}}]
  \label{thm:intro-tilting-hypersurface}
  Let $A$ be a noncommutative isolated quadric hypersurface singularity
  of dimension $d\geq 2$.
  Let $\{M_1, \dotsc, M_n\}$ be a complete set of non-isomorphic
  indecomposable non-projective graded Cohen--Macaulay $A$-modules
  generated in degree~0. Then
  \(
    (\pi(A)(-d+2),\pi(A)(-d+3),\dotsc,\pi(A),\pi(M_1),\dotsc,
    \pi(M_n))
  \)
  is a full strong exceptional sequence in $\Dbqgr{A}$.
  
  Consequently, there is a triangulated equivalence 
  \(
    \Dbqgr{A}\simeq \Dbmod{\Lambda},
  \)
  where $\Lambda \defeq
  \End_{\gr{A}}\bigl(A(-d+2)\oplus\dotsb\oplus A(-1)\oplus A
  \oplus M_1\oplus\dotsb\oplus M_n\bigr)$ is the endomorphism algebra, and $\Lambda$ 
  is finite-dimensional and basic.
\end{theorem}

In \ref{sec:ct-hyper}, we study the existence of graded $(d-1)$-cluster tilting modules
over noncommutative isolated quadric hypersurface singularities $S/(f)$ of dimension $d\geq 2$.
If $f$ is central, no such module exists for $d\geq 4$
(\ref{cor:central-no-NCCR}),
and therefore no noncommutative resolution exists;
this includes the case when $A=S/(f)$ is commutative (\ref{rem:comm-no-NCCR}).
For noncentral normal regular $f$, we construct examples admitting such modules
in infinitely many dimensions.
The following theorem characterizes precisely which dimensions can occur.

\begin{theorem}[{\ref{thm:ct-dimension}}]
  \label{thm:intro-ct-dimension}
  There exists
  a noncommutative isolated quadric hypersurface singularity
  of dimension $d\geq 2$ with a graded $(d-1)$-cluster tilting module if and only if
  $d=2^n+1$ for some integer $n\geq 0$.
\end{theorem}


The proof of \ref{thm:intro-ct-dimension} relies on the following criterion.
Let $A$ be a noncommutative isolated quadric hypersurface singularity
of dimension $d\geq 2$.
The functor $\rho\defeq\Omega(\blank)(1)$, where $\Omega$ denotes
the syzygy functor, induces a permutation of the finite set
of isomorphism classes of indecomposable non-projective graded
Cohen--Macaulay $A$-modules generated in degree 0.

\begin{proposition}[{\ref{prop:ct-decomp}}]
  The algebra $A$ admits a graded $(d-1)$-cluster tilting module
  if and only if the cardinality of every $\rho$-orbit
  is divisible by $d-1$.
\end{proposition}

If $A$ admits such a graded $(d-1)$-cluster tilting module,
then $d-1$ divides $2^d$, forcing $d=2^n+1$ for some $n\geq 0$.
For $d=2$, the homogeneous coordinate ring
of a smooth projective conic provides an example.
For each $d=2^n+1$ with $n\geq 1$, we construct a graded twist
of a $(-1)$-skew quadric hypersurface singularity for which every $\rho$-orbit
has cardinality $d-1$ (\ref{ex:ct}).

In \ref{sec:ct-surface}, we consider noncommutative isolated quadric surface singularities.
In a series of papers
\cite{huNoncommutativeConicsCalabi2023,
huNoncommutativeAffinePencils2026,
huClassificationNoncommutativeCentral2026},
the authors obtained a complete classification of
noncommutative projective central conics (the case $d=2$) in characteristic 0.
Their classification relied on the classification of 3-dimensional quantum polynomial algebras $S$ and explicit descriptions of the central quadratic elements $f$.
Classifying the homogeneous coordinate algebras up to graded isomorphism is a finer problem,
since two non-isomorphic noncommutative isolated quadric surface singularities
may define isomorphic noncommutative projective conics.
We classify noncommutative isolated quadric surface singularities, including those defined by normal quadratic elements, up to graded isomorphism over an algebraically closed field of characteristic $\neq$ 2.
We do so by classifying the associated noncommutative resolutions
and recovering the original graded algebras from them,
thereby avoiding a case-by-case analysis of the pair $(S,f)$.

Let $A$ be a noncommutative isolated quadric surface singularity,
and let $M_1,\dotsc, M_n$ and $\Lambda$ be as in \ref{thm:intro-tilting-hypersurface}.
Then $A\oplus M_1\oplus\dotsb\oplus M_n$ is a graded 1-cluster tilting module.
Its graded endomorphism algebra is a noncommutative resolution of $A$
whose degree-zero part is $\Lambda$.
By \ref{enum:prop-NCCR-iso}, each such resolution is a graded twist
of the graded preprojective algebra of a quiver $Q$ (\ref{dfn:preproj-alg}),
where $Q$ is a Euclidean quiver
of type $\widetilde{A_1}$ or $\widetilde{D_4}$
with the orientation specified in \ref{ex:A1,ex:D4}.
These algebras are called twisted preprojective algebras of type $Q$
(\ref{dfn:twisted-preproj}),
which are graded skew Calabi--Yau of dimension 2 with Gorenstein parameter 1 (\ref{cor:twisted-preproj}).

Using the Morita-like theory for graded quotient categories $\qgrOp$
developed in \cite{liNoncommutativeResolutionsASGorenstein2026}
and further studied in \cite{liNoncommutativeResolutionsNoncommutative2026},
we establish the following correspondence.

\begin{theorem}[{\ref{thm:classif-corresp}}]
  \label{thm:intro-classif-corresp}
  There is a one-to-one correspondence
  between the graded isomorphism classes of
  noncommutative isolated quadric surface singularities
  and the graded isomorphism classes of
  twisted preprojective algebras of type $\widetilde{A_1}$ or $\widetilde{D_4}$.
\end{theorem}

By classifying the twisted preprojective algebras (\ref{thm:B-classif})
and using \ref{thm:intro-classif-corresp}, we obtain the following.

\begin{theorem}[{\ref{thm:A-classif}}]
  \label{thm:intro-A-classif}
  There is a one-to-one correspondence between the graded isomorphism classes of
  noncommutative isolated quadric surface singularities
  and the union of the sets of conjugacy classes in $\PGL_2(\kk)$ and $S_4$.
\end{theorem}

In fact, \ref{thm:intro-A-classif} gives a classification of noncommutative isolated quadric surface singularities,
which consists of one family and six isolated cases; see \ref{tab:A-classif}.

We also determine which algebras in \ref{tab:A-classif}
admit a presentation by a central quadratic element
(\ref{prop:A-classif-central}).
Exactly two cases, $\mathbb{D}_4$ and $\mathbb{D}_5$,
do not admit such a presentation.
All remaining cases are graded isomorphic to noncommutative central quadric surface singularities,
even when the defining element is not central.



\section{Preliminaries}\label{sec:prelim}

Throughout this paper, we work over an algebraically closed field $\kk$ with $\charOp(\kk)\neq 2$.
By an algebra, we mean a $\kk$-algebra; by an $A$-module, we mean a right $A$-module.

For a $\ZZ$-graded algebra $A$, let $\Gr{A}$ denote the category of
$\ZZ$-graded $A$-modules with degree-zero morphisms,
and let $\gr{A}$ denote its full subcategory of finitely generated modules.
For any $M\in\Gr{A}$, the $j$-th shift $M(j)$ is defined by
$M(j)_i=M_{i+j}$. For $M, N\in\Gr{A}$, define
\[
  \gExt_A^n(M,N)=\bigoplus_{j\in\ZZ}\Ext_{\Gr{A}}^n(M,N(j)).
\]
For $M\in\gr{A}$, let $\add(M)$ denote the full subcategory
of $\gr{A}$ consisting of direct summands of finite direct sums
of degree shifts of $M$.


From now on,
every $\ZZ$-graded algebra $A$ is assumed to be locally finite,
that is, $\dim_{\kk}(A_i)<\infty$ for every $i\in\ZZ$.
Let $A$ be a noetherian $\ZZ$-graded algebra that is bounded below,
and let $J_A$ be its graded Jacobson radical.
For $M\in\Gr{A}$, the torsion functor $\Gamma_A$ is given by
\[
  \Gamma_A(M)=\{m\in M\mid mJ_A^n=0 \text{ for some }n\geq 1\}.
\]
If $\Gamma_A(M)=M$, then $M$ is called torsion.
Let $\Tors{A}$ be the full subcategory of $\Gr{A}$ consisting of torsion $A$-modules,
and let $\tors{A}=\Tors{A}\cap\gr{A}$. Then $\tors{A}$ coincides with the full subcategory of finite-dimensional graded $A$-modules.
Let $\QGr{A}=\Gr{A}/\Tors{A}$ and let $\qgr{A}=\gr{A}/\tors{A}$.
Both the quotient functors are denoted by $\pi$.

The right derived functors $\uR^i\Gamma_A$ of $\Gamma_A$ are called the local cohomology functors of $A$.
We define
\[
  \dep_A(M)
  =\inf\{i\geq 0\mid \gExt_A^i(A/J_A,M)\neq 0\},
\]
which equals $\inf\{i\geq 0\mid \uR^i\Gamma_A(M)\neq 0\}$
by \cite[Lem.\@ 2.12]{liCommonlyGradedAlgebras2025}.
Suppose that $\Gamma_A$ has finite cohomological dimension $d$, that is,
\[d=\sup\{i\geq 0\mid\uR^i\Gamma_A(M)\neq 0\text{ for some }M\in\Gr{A}\}<\infty.\]
A nonzero module $M\in\gr{A}$ is called \dfn{Cohen--Macaulay}
if $\dep_A(M)=d$.
We denote by $\CM{A}$ the full subcategory consisting of Cohen--Macaulay $A$-modules
and the zero module, and by $\CM*{A}$ its full subcategory consisting
of modules generated in degree 0.
The \dfn{stable category} $\sCM{A}$ has the same objects as $\CM{A}$, with
morphisms given by
\[
  \Hom_{\sCM{A}}(M, N)\defeq \Hom_{\gr{A}}(M,N) / P(M, N),
\]
where $P(M, N)$ consists of the degree-zero $A$-module maps $f\colon M\to N$
that factor through a projective object in $\Gr{A}$.
We denote by $\sCM*{A}$ the full subcategory of $\sCM{A}$ consisting of modules generated
in degree 0.

By \cite[Lem.\@ 3.4]{smithNoncommutativeQuadricSurfaces2013},
two modules $M,N\in\CM{A}$ with no nonzero projective direct summands
are isomorphic in $\sCM{A}$ if and only if they are isomorphic in $\gr{A}$.
For an additive category $\cC$, we denote by $\ind(\cC)$ the set
of isomorphism classes of indecomposable objects.



\subsection{Some properties of graded algebras}

We first recall the definitions of Koszul algebras, Artin--Schelter Gorenstein and regular algebras, and $\NN$-graded skew Calabi--Yau algebras as a foundation for studying noncommutative quadric hypersurface singularities and their noncommutative resolutions.

Let $A$ be a connected graded algebra, that is, an $\NN$-graded algebra with $A_0=\kk$.
A graded $A$-module $M$ is said to have a \dfn{linear resolution} if it has a graded free resolution
$\dotsb\to P_2\to P_1\to P_0\to M\to 0$
such that $P_i$ is generated in degree $i$ for every $i\geq 0$. The algebra
$A$ is called \dfn{Koszul} if $\kk_A$ has a linear resolution.
If $A$ is a Koszul algebra, then
$A^!\defeq\bigoplus_{i\in\NN}\gExt_A^i(\kk,\kk)$
is also a Koszul algebra and is called the \dfn{Koszul dual} of $A$.
We refer the reader to \cite{beilinsonKoszulDualityPatterns1996} for details.

\begin{definition}
  \label{dfn:AS-Gorenstein}
  A connected graded algebra $A$ is called \dfn{Artin--Schelter Gorenstein} (or \dfn{AS-Gorenstein}, for short) of dimension $d$ if
  $\idim_A(A)=\idim_{\Ao}(A)=d<\infty$, and, for some integer $l$,
  \[
    \gExt_A^i(\kk,A)\cong\gExt_{\Ao}^i(\kk,A)\cong
    \begin{cases}
      \kk(l),&i=d,\\
      0,&i\neq d.
    \end{cases}
  \]
  The integer $l$ is called the \dfn{Gorenstein parameter} of $A$.
  If, in addition, $\gldim(A)=d$, then $A$ is called \dfn{Artin--Schelter regular} (or \dfn{AS-regular}, for short) of dimension $d$.
\end{definition}

A \dfn{quantum polynomial algebra of dimension $d$} is
a noetherian AS-regular algebra $S$
of dimension $d$ with Hilbert series $H_S(t)=(1-t)^{-d}$.
It is Koszul (\cite[Thm.\@ 5.11]{smithFinitedimensionalAlgebrasRelated1996}),
so its Gorenstein parameter is $d$.




The noncommutative resolutions considered in this paper are not necessarily connected graded.
We therefore recall the notion of AS-regularity over the degree-zero part of an $\NN$-graded algebra.

Let $A$ be a $\ZZ$-graded algebra and let $\sigma$ be a graded automorphism of $A$.
For a graded $A$-module $M$,
we write $M^\sigma$ for the graded right $A$-module
with the same underlying graded $\kk$-vector space as $M$
and right action defined by
$m\cdot^\sigma a=m\sigma(a)$
for all $m\in M$ and $a\in A$.
Let $D\defeq\gHom_{\kk}(\blank,\kk)$ be the Matlis duality.

\begin{definition}[{\cite[Def.\@ 3.1]{minamotoStructureASGorensteinAlgebras2011}}]
  \label{dfn:AS-regular-over-zero}
  Let $A$ be an $\NN$-graded algebra and let $R=A_0$.
  The algebra $A$ is called \dfn{AS-regular over $R$ of dimension $d$
  with Gorenstein parameter $l$} if $\gldim(A)=d$, and
  \[
    \RgHom_A(R,A)\cong{}^1D(R)^\nu(l)[-d]
  \]
  in the derived category of graded $(R,R)$-bimodules
  for some $\nu\in\Aut{R}$.
\end{definition}

The original definition also requires $\gldim(R)<\infty$,
but this condition is redundant by \cite[Lem.\@ 3.23]{reyesGradedTwistedCalabi2022}.


\begin{definition}[{{\cite[Def.\@ 1.2]{reyesGradedTwistedCalabi2022}}}]
\label{dfn:CY}
  Let $A$ be an $\NN$-graded algebra.
  \begin{enumerate}
    \item $A$ is \dfn{graded homologically smooth} if it is graded perfect
      as an $\Ae$-module, i.e., it has a finitely generated
      graded $\Ae$-projective resolution of finite length.
    \item $A$ is \dfn{graded twisted Calabi--Yau (CY) of dimension $d$}
      if it is graded homologically smooth and
      $\RgHom_{\Ae}(A,\Ae)\cong U[-d]$
      in the derived category of graded $(A,A)$-bimodules,
      for some graded invertible bimodule $U$.
    \item $A$ is \dfn{graded skew CY}
      if it is graded twisted CY and
      $U\cong {}^1A^{\mu}(l)$ as graded $(A,A)$-bimodules
      for some $\mu\in\GrAut{A}$ and some integer $l$. In this case,
      $\mu$ is called the \dfn{Nakayama automorphism} of $A$, and
      $l$ is called the \dfn{Gorenstein parameter} of $A$.
    \item $A$ is \dfn{graded CY}
      if it is graded skew CY and $\mu$ is inner.
  \end{enumerate}
\end{definition}

\begin{theorem}[{{\cite[Thms.\@ 5.17, 6.4]{liCommonlyGradedAlgebras2025}}}]
  \label{thm:CY-equiv}
  Let $A$ be an $\NN$-graded algebra such that $A_0$ is basic.
  The following are equivalent.
  \begin{enumerate}
    \item $A$ is graded skew CY of dimension $d$
      with Gorenstein parameter $l$.
    \item $A$ is $\NN$-graded AS-regular of dimension $d$
      with constant Gorenstein parameters $(l, l, \dotsc, l)$
      in the sense of \cite{liCommonlyGradedAlgebras2025}. 
    \item $A$ is AS-regular over $A_0$ of dimension $d$
      with Gorenstein parameter~$l$.
  \end{enumerate}
\end{theorem}

\begin{remark}
  Different definitions of the Nakayama automorphism appear in the literature.
  We follow the convention of
  \cite[Def.\@ 4.4]{brownDualisingComplexesTwisted2008},
  \cite[Def.\@ 0.1(b)]{reyesSkewCalabiYau2014},
  \cite[Def.\@ 1.2(iii)]{reyesGradedTwistedCalabi2022},
  and \cite[Defs.\@ 5.7 and 6.1]{liCommonlyGradedAlgebras2025}.
  The Nakayama automorphism in
  \cite[Def.\@ 3.9]{minamotoStructureASGorensteinAlgebras2011}
  and \cite[Sec.\@ 2.2]{ueyamaDerivedCategoriesSkew2023}
  corresponds to the inverse of what we define here.
  See \cite[Prop.\@ 4.3]{brownDualisingComplexesTwisted2008}
  and \cite[Prop.\@ 4.14]{liCommonlyGradedAlgebras2025} for details.
\end{remark}

\subsection{Twists of graded algebras}

The twist of graded algebras is an important construction,
first investigated in \cite[Sec.\@ 8]{artinModulesRegularAlgebras1991}
and subsequently extended in \cite{zhangTwistedGradedAlgebras1996}.
We recall this construction and its basic properties
for use in our study of cluster tilting modules and noncommutative resolutions.
For simplicity, we consider only the $\NN$-graded case,
although the construction applies to any algebra graded by a semigroup.

\begin{definition}
  A \dfn{twisting system} on an $\NN$-graded algebra $A$
  is a sequence
  $\theta=\{\theta_i\}_{i\in\ZZ}$
  of graded $\kk$-linear bijections of $A$ such that
  \begin{equation*}
    \theta_i(a\theta_j(b))=\theta_i(a)\theta_{i+j}(b)
    \quad\text{for all }a\in A_j,b\in A \text{ and }
    i,j\in\ZZ.
  \end{equation*}
  The \dfn{twisted algebra} $A_{\theta}$ has the same underlying graded $\kk$-vector space as $A$. Its elements are denoted by $a_\theta$ for $a\in A$,
  and multiplication is defined by
  \begin{equation*}
    a_\theta b_\theta\defeq (a\theta_{i}(b))_\theta
    \quad\text{for }a\in A_i,b\in A.
  \end{equation*}

  For $M\in\Gr{A}$, the \dfn{twisted module} $M_\theta$
  has the same underlying graded $\kk$-vector space as $M$,
  with elements denoted by $m_\theta$ for $m\in M$,
  and an $A_\theta$-action defined by
  \[
    m_\theta a_\theta\defeq (m\theta_i(a))_\theta
    \quad\text{for }m\in M_i,a\in A.
  \]
  \end{definition}

  If $\tau$ is a graded automorphism of $A$,
  then $\theta\defeq\{\tau^i\}_{i\in\ZZ}$
  is a twisting system on $A$.
  In this case, the twisted algebra $A_{\theta}$
  is denoted by $A_{\tau}$.

By \cite[Prop.\@ 2.5]{zhangTwistedGradedAlgebras1996},
if $A$ is an $\NN$-graded algebra and $\theta$ is a twisting system on $A$,
then $\theta^{-1}\defeq\{\theta_i^{-1}\}_{i\in\ZZ}$ is a twisting system
on $B\defeq A_\theta$ such that $A=B_{\theta^{-1}}$.

Let $S$ be an $\NN$-graded algebra,
and let $F=\{f_\lambda\}_{\lambda\in\Lambda}$ be a set of homogeneous elements of $S$.
Let $I$ and $V$ denote, respectively, the two-sided ideal and the vector space generated by $F$,
and let $A=S/I$ be the quotient algebra with quotient map $\pi\colon S\to A$.
Suppose that $\Tau$ is a graded automorphism of $S$ such that $\Tau(V)=V$.
Then $\Tau$ induces an automorphism $\tau$ of $A$
such that $\tau\pi=\pi\Tau$.

The following lemma describes the twist of the quotient algebra $S/I$.
It is a generalization of
\cite[Lem.\@ 4.16]{huClassificationNoncommutativeCentral2026}.
\begin{lemma}\label{lem:twisted-alg}
  Keep the notation as above. Then $A_\tau=S_\Tau/(F_\Tau)$,
  where $F_\Tau$ denotes the subset $\{g_\Tau\mid g\in F\}$ of $S_\Tau$.
\end{lemma}

\begin{proof}
  By definition, there is a natural projection
  $\pi\colon S_\Tau\to A_\tau$ with $\ker(\pi)=I_\Tau$.
  Thus, $I_\Tau$ is a two-sided ideal of $S_\Tau$,
  and it suffices to show $F_\Tau$ generates $I_\Tau$ as a two-sided ideal.
  For any $f_\Tau\in I_\Tau$,
  write $f$ as a finite sum $f=\sum_i x_ig_iy_i$,
  where $g_i\in F$ and $x_i, y_i\in S$ are homogeneous.
  Then
  $f_\Tau=\sum_i (x_i)_\Tau (\Tau^{-\abs{x_i}}(g_i))_\Tau
  (\Tau^{-\abs{x_ig_i}}(y_i))_\Tau\in (V_\Tau)=(F_\Tau)$
  as $\Tau(V)=V$.
\end{proof}

We record several properties of graded algebras and modules that are preserved under twisting in the following propositions.

\begin{proposition}[{{\cite{zhangTwistedGradedAlgebras1996}}}]
  \label{prop:zhang-twist-1}
  Let $A$ be an $\NN$-graded algebra, let
  $\theta=\{\theta_i\}_{i\in\ZZ}$ be a twisting system on $A$,
  and let $B=A_\theta$.
  The following hold.
  \begin{enumerate}[ref={\ref{prop:zhang-twist-1}(\arabic*)}]
    \item\label{enum:prop-zhang-twist-equiv}
      \cite[Thms.\@ 3.1 and 3.4]{zhangTwistedGradedAlgebras1996}
      The assignment $F(M)=M_\theta$, with $F$ acting on morphisms
      via the same underlying linear maps,
      defines an equivalence
      $F\colon\Gr{A}\simeq\Gr{B}$.
      For every $j\in\ZZ$, there is a graded $B$-module isomorphism
      \[
        F(A(j))\cong B(j),
        a_\theta\mapsto(\theta_{-j}^{-1}(a))_\theta.
      \]
    \item\label{enum:prop-zhang-twist-noeth}
      \cite[Thm.\@ 3.7(1) and Prop.\@ 5.1]{zhangTwistedGradedAlgebras1996}
      $A$ is noetherian if and only if $B$ is.
      In this case, $F$ restricts to an equivalence $\gr{A}\simeq\gr{B}$
      and induces an equivalence $\qgr{A}\simeq\qgr{B}$.
    \item\label{enum:prop-zhang-twist-AS}
      \cite[Thm.\@ 5.11]{zhangTwistedGradedAlgebras1996}
      If $A$ is noetherian, then
      $A$ is AS-Gorenstein (resp.\@ AS-regular) of dimension $d$
      with Gorenstein parameter $l$ if and only if $B$ is.
  \end{enumerate}
\end{proposition}

\begin{proposition}\label{prop:zhang-twist-2}
  Keep the notation in \ref{prop:zhang-twist-1}.
  The following hold.
  \begin{enumerate}[ref={\ref{prop:zhang-twist-2}(\arabic*)}]
    \item\label{enum:prop-zhang-twist-CY} $A$ is graded basic if and only if $B$ is.
      If $A$ is graded basic, then $A$ is graded skew CY
      of dimension $d$ with Gorenstein parameter $l$
      if and only if $B$ is.
    \item\label{enum:prop-zhang-twist-Koszul} If $A$ is connected graded, then $A$ is Koszul
      if and only if $B$ is.
      Moreover, $A$ is a quantum polynomial algebra of dimension $d$
      if and only if $B$ is.
    \item\label{enum:prop-zhang-twist-CM} If $A$ is noetherian AS-Gorenstein, then the functor $F$ restricts to an equivalence
      $\CM{A}\simeq\CM{B}$, and induces equivalences
      $\sCM{A}\simeq\sCM{B}$, $\sCM*{A}\simeq\sCM*{B}$.
  \end{enumerate}
\end{proposition}

\begin{proof}
  Since $A\cong B_{\theta^{-1}}$, it suffices to prove the implications
  from $A$ to $B$.

  (1)
  By \cite[Prop.\@ 2.4]{zhangTwistedGradedAlgebras1996},
  we may assume without loss of generality that $\theta_0=\id_A$.
  Then $A_0=B_0$,
  so $A$ is graded basic if and only if $B$ is.

  Suppose that $A$ is graded basic and graded skew CY
  of dimension $d$ with Gorenstein parameter $l$.
  Let $R=A_0=B_0$ and write $D(R)=\Hom_{\kk}(R,\kk)$.
  By \ref{thm:CY-equiv},
  \[
    \RgHom_A(R,A)\cong{}^1D(R)^\nu(l)[-d]
  \]
  for some $\nu\in\Aut{R}$.
  Let $i\geq0$ and $j\in\ZZ$.
  Since $\theta_0=\id_A$,
  $\theta_{-j}$ restricts to an algebra automorphism of $R$.
  By \ref{enum:prop-zhang-twist-equiv},
  there is an isomorphism of $(R,R)$-bimodules
  \[
    {}^{\theta_{-j}\vert_R}\Ext^i_{\Gr{A}}(R,A(j))
    \cong\Ext^i_{\Gr{B}}(R,B(j)).
  \]
  The left $R$-actions on the Ext groups are induced by left
  multiplication on $A(j)$ and $B(j)$, respectively.
  Let $\alpha=\theta_l\vert_R$. Then
  \[
    \gExt_B^i(R,B)\cong
    \begin{cases}
      {}^\alpha D(R)^\nu(l)
      \cong{}^1D(R)^{\alpha^{-1}\nu}(l),&i=d,\\
      0,&i\neq d.
    \end{cases}
  \]
  Thus, $B$ is AS-regular over $R$ of dimension $d$
  with Gorenstein parameter $l$, and \ref{thm:CY-equiv} applies.

  (2) If $A$ is Koszul, then $F$ sends a linear resolution of $\kk_A$
  to a linear resolution of $\kk_B$,
  as $F(\kk_A)=\kk_B$ and $F(A(j))\cong B(j)$
  by \ref{enum:prop-zhang-twist-equiv}.
  The assertion on quantum polynomial algebras follows from
  (1), \ref{enum:prop-zhang-twist-noeth}, and $H_A(t)=H_B(t)$.

  (3) Suppose $M\in\gr{A}$, $i\geq0$, and $j\in\ZZ$.
  It follows from $F(\kk_A(j))=\kk_B(j)$ that
  \[
    \Ext^i_{\Gr{A}}(\kk_A(j),M)
    \cong\Ext^i_{\Gr{B}}(\kk_B(j),F(M)).
  \]
  Hence $\dep_A(M)=\dep_B(F(M))$, so $F$ restricts
  to an equivalence $\CM{A}\simeq\CM{B}$
  by \ref{enum:prop-zhang-twist-AS}.
  Since $F$ preserves projective modules and modules generated
  in degree 0, the remaining assertions follow.
\end{proof}

We also need to determine when a normal regular element
remains normal and regular after twisting.

\begin{lemma}\label{lem:twist-normal-regular}
  Let $A$ be an $\NN$-graded algebra, let $f\in A$ be homogeneous of degree $s$,
  and let $\tau$ be a graded automorphism of $A$.
  \begin{enumerate}[beginpenalty=10000,midpenalty=10000]
    \item If $f$ is regular in $A$, then $f_\tau$ is regular in $A_\tau$.
    \item If $f$ is normal in $A$ with $\tau(f)\in \kk f$,
      then $f_\tau$ is normal in $A_\tau$.
  \end{enumerate}
\end{lemma}

\begin{proof}
  (1) is immediate.
  For (2), a direct computation shows that
  \begin{align*}
    (A_\tau)_if_\tau&=(A_i\tau^i(f))_\tau=(A_if)_\tau;\\
    f_\tau(A_\tau)_i&=(f\tau^s(A_i))_\tau=(fA_i)_\tau=(A_if)_\tau.
  \end{align*}
  Hence, $A_\tau f_\tau=f_\tau A_\tau$.
\end{proof}

\subsection{Noncommutative quadric hypersurface singularities}


We recall the basic properties of
noncommutative quadric hypersurface singularities and
their Cohen--Macaulay modules.

\begin{definition}\label{dfn:quadric}
  A \dfn{noncommutative quadric hypersurface singularity of dimension $d$}
  is a graded algebra $A=S\div (f)$,
  where $S$ is a quantum polynomial algebra of dimension $d+1$
  and $f\in S_2$ is a normal regular element;
  $A$ is called a \dfn{noncommutative quadric surface singularity}
  if $d=2$.
  The noncommutative quadric (hyper)surface singularity $A=S\div (f)$ is called \dfn{central} if $f\in S$ is central.
\end{definition}

The normal regular element $f \in S$ induces a graded automorphism $\sigma$ of $S$ such that $fa=\sigma(a)f$ for all $a\in S$. This automorphism is called the \dfn{normalizing automorphism} induced by $f$.
Since $\sigma(f)=f$, $\sigma$ also induces a graded automorphism of $A$, which is again denoted by $\sigma$.

Any noncommutative quadric hypersurface singularity $A$ of dimension $d$
is AS-Gorenstein of dimension $d$ with Gorenstein parameter $l=d-1$, by Rees' Lemma (see, for example,
\cite[Prop.\@ 3.4(b)]{levasseurPropertiesNoncommutativeRegular1992}).
Moreover, since $l<d$, $A$ has infinite global dimension by \cite[Prop.\@ 3.1(4)]{stephensonGrowthGradedNoetherian1997}.

\begin{lemma}\label{lem:normal-quadric-dual}
  Let $A=S\div(f)$ be a noncommutative quadric hypersurface singularity of dimension $d$. Then
  \begin{enumerate}[ref={\ref{lem:normal-quadric-dual}(\arabic*)},
    beginpenalty=10000,midpenalty=10000]
    \item\label{enum:hilbert-quadric}
      $A$ has Hilbert series $H_A(t)=(1-t)^{-d}(1+t)$.
    \item $A$ is Koszul.
    \item\label{enum:w}
      There exists a normal regular element
      $w\in A_2^!$ such that $A^!\div(w)\cong S^!$.
    \item\label{enum:w-central}
      If $f\in S$ is central, then the element $w\in A^!$ in (3)
      can be chosen to be central.
  \end{enumerate}
\end{lemma}

\begin{proof}
  The exact sequence $0\to S(-2)\to S\to A\to 0$
  of graded $S$-modules implies that $H_A(t)=(1-t^2)H_S(t)=(1-t)^{-d}(1+t)$, proving (1).
  Statements (2) and (3) follow from \cite[Thm.\@ 1.2]{sheltonKoszulAlgebrasNew2001}
  and \cite[Cor.\@ 1.4(1)]{sheltonKoszulAlgebrasNew2001}, respectively.
  Statement (4) is asserted in
  \cite[Lem.\@ 5.1(2)]{smithNoncommutativeQuadricSurfaces2013},
  and a detailed proof is given in
  \cite[pp.\@ 72--73]{heCliffordDeformationsKoszul2024}.
\end{proof}

The graded Cohen--Macaulay representations of
a noncommutative quadric hypersurface singularity $A$
are controlled by the finite-dimensional algebra
\begin{equation*}
  C(A)\defeq A^![w^{-1}]_0,
\end{equation*}
the degree-zero part of the localization of $A^!$
with respect to the element $w\in A_2^!$ given in \ref{enum:w},
as summarized in the following lemma.

\begin{lemma}
  [{{\cite[Lem.\@ 5.1, Prop.\@ 5.2]{smithNoncommutativeQuadricSurfaces2013};
    \cite[Lem.\@ 4.13]{moriNoncommutativeKnorrersPeriodicity2022}}}]
  \label{lem:CA}
  Let $A$ be a noncommutative quadric hypersurface singularity
  of dimension $d$. Then
  \begin{enumerate}[ref={\ref{lem:CA}(\arabic*)}]
    \item\label{enum:lem-CA-dim} $\dim_{\kk}(C(A))=2^d$.
    \item\label{enum:lem-CA-F} There is an equivalence of triangulated categories
      $F\colon \sCM{A}\simeq \Dbmod{C(A)}$.
    \item\label{enum:lem-CA-G} $F$ restricts to an equivalence
      $G\colon \sCM*{A}\simeq \mod{C(A)}$.
  \end{enumerate}
\end{lemma}

Recall that the global dimension of an abelian category $\cC$
is defined by
\begin{equation*}
  \gldim(\cC)\defeq\sup\{n\mid\Ext^n_{\cC}(\cM, \cN)\neq 0
  \text{ for some }\cM,\cN\in\cC\}.
\end{equation*}

Noncommutative graded isolated singularities were first defined in
\cite{ueyamaGradedMaximalCohen2013}.
\begin{definition}\label{dfn:iso-sing}
  Let $A$ be a noetherian AS-Gorenstein algebra
  of dimension $d\geq 1$. If $\gldim(\qgr{A})<\infty$, then
  $A$ is called a \dfn{graded isolated singularity}.
\end{definition}

Suppose that $A$ is a graded isolated singularity.
By \cite[Prop.\@ A.3]{reitenNoetherianHereditaryAbelian2002}, 
\cite[Cor.\@ 4.5]{ueyamaGradedMaximalCohen2013}, and 
\cite[Cor.\@ 3.9]{iyamaCohenMacaulayRepresentationsArtinSchelter2025}, 
\begin{equation}\label{eq:iso-sing}
  \gldim(\qgr{A})=\gldim(\QGr{A})=d-1.
\end{equation}
If $\add(M)=\CM{A}$ for some $M\in\CM{A}$, then $A$ is said to have \dfn{finite graded Cohen--Macaulay representation type}.
Such an $M$ is called an \dfn{additive generator} of $\CM{A}$.

\begin{theorem}
  [{{\cite[Thm.\@ 5.6]{smithNoncommutativeQuadricSurfaces2013};
    \cite[Thm.\@ 5.5]{moriNoncommutativeKnorrersPeriodicity2022}}}]
  \label{thm:quadric-iso-sing}
  Let $A$ be a noncommutative quadric hypersurface singularity of dimension $d\geq 1$.
  Then the following are equivalent.
  \begin{enumerate}
    \item $A$ is a graded isolated singularity.
    \item $A$ has finite graded Cohen--Macaulay representation type.
    \item $C(A)$ is semisimple.
  \end{enumerate}
\end{theorem}

\begin{lemma}
  [{{\cite[Lem.\@ 5.2]{moriNoncommutativeKnorrersPeriodicity2022}}}]
  \label{lem:CM-shift}
  Let $A$ be a noncommutative isolated quadric hypersurface singularity.
  Then every indecomposable graded Cohen--Macaulay
  $A$-module is isomorphic to $M(t)$ for some
  $M\in\ind(\CM*{A})$ and some $t\in\ZZ$.
\end{lemma}


\section{Tilting for
noncommutative isolated quadric hypersurface singularities}\label{sec:tilting}

If $A$ is a commutative isolated quadric hypersurface singularity, then
$X=\Proj{A}$ is smooth and $\Dbqgr{A}\simeq\Dbcoh{X}$.
Kapranov showed that $\Dbcoh{X}$ admits a full strong exceptional sequence
(see \cite{kapranovDerivedCategoryCoherent1986} and \cite[Sec.\@ 4]{kapranovDerivedCategoriesCoherent1988}).
This sequence consists of twists of the structure sheaf and spinor bundles.
The spinor bundles are sheaves associated to graded Cohen--Macaulay $A$-modules.


In the noncommutative setting,
Ueyama \cite[Thm.\@ 3.21]{ueyamaDerivedCategoriesSkew2023}
generalized this construction to
$(\pm 1)$-skew quadric hypersurface singularities.
Iyama--Kimura--Ueyama
\cite[Thm.\@ 9.8(3)]{iyamaCohenMacaulayRepresentationsArtinSchelter2025}
constructed a full strong exceptional sequence
for arbitrary noncommutative isolated quadric hypersurface singularities
using shifted syzygies of the trivial module.

For the application to noncommutative resolutions in \ref{sec:ct-surface},
we construct a tilting object of the form $\pi(M)$ in $\Dbqgr{A}$,
where $M$ is a Cohen--Macaulay generator of $A$.
In \ref{thm:tilting-hypersurface}, we achieve this
by generalizing Kapranov's and Ueyama's exceptional sequences
to arbitrary noncommutative isolated quadric hypersurface singularities.

\begin{definition}
  Let $\cT$ be a $\kk$-linear triangulated category, and let $E\in\cT$.
  The object $E$ is called \dfn{exceptional} if $\Hom(E, E[p])=0$ for $p\neq 0$ and $\Hom(E,E)=\kk$.
  An \dfn{exceptional sequence} in $\cT$ is a sequence of exceptional objects
  $(E_0, \dotsc, E_n)$ satisfying the semi-orthogonality condition
  $\Hom(E_i, E_j[p])=0$ for all $p$ whenever $i>j$.
  An exceptional sequence is called \dfn{strong}
  if $\Hom(E_i, E_j[p])=0$ for all $i$ and $j$ whenever $p\neq 0$.
  An exceptional sequence is called \dfn{full}
  if it generates $\cT$ as a triangulated category.
\end{definition}

The following form of Serre duality is useful for computing $\Ext$ groups in $\qgr{A}$.
Recall that a noetherian AS-Gorenstein algebra $A$ of dimension $d$ with Gorenstein parameter $l$
has a balanced dualizing complex
${}^1A^{\mu^{-1}}(-l)[d]$ for some $\mu\in\GrAut{A}$
(\cite{yekDualizingComplexesNoncommutative1992}).
The automorphism $\mu$ is called the \dfn{Nakayama automorphism} of $A$.

\begin{lemma}[{{\cite[Lem.\@ 2.3]{ueyamaDerivedCategoriesSkew2023}}}]
  \label{lem:Serre-duality}
  Let $A$ be an AS-Gorenstein isolated singularity
  of dimension $d\geq 2$ with Gorenstein parameter $l$
  and Nakayama automorphism $\mu$.
  Then, for $M,N\in\CM{A}$, the following statements hold.
  \begin{enumerate}
    \item $\Ext_{\qgr{A}}^q(\pi(M),\pi(N))\cong\Ext_{\gr{A}}^q(M,N)$
      for $0\leq q\leq d-2$.
    \item $\Ext_{\qgr{A}}^{d-1}(\pi(M),\pi(N))
      \cong D\bigl(\Hom_{\gr{A}}(N, M^{\mu^{-1}}(-l))\bigr)$.
    \item $\Ext_{\qgr{A}}^q(\pi(M),\pi(N))=0$ for $q\geq d$.
  \end{enumerate}
\end{lemma}

Constructing an exceptional sequence over a noncommutative quadric hypersurface singularity requires the following properties of Cohen--Macaulay modules generated in degree 0. These properties are also used in \ref{sec:ct-hyper}.
\begin{lemma}\label{lem:CM-module-quadric}
  Let $A=S\div (f)$ be a noncommutative quadric hypersurface singularity,
  and let $0\neq M\in\CM*{A}$ have no nonzero projective direct summands.
  Set $s=\dim_{\kk}(M_0)$. Then
  \begin{enumerate}[beginpenalty=10000,midpenalty=10000,
    ref={\ref{lem:CM-module-quadric}(\arabic*)}]
    \item 
      $\Omega^2(M)\cong M^{\sigma}(-2)$, where $\sigma$ is the normalizing automorphism induced by $f$.
    \item\label{enum:lem-CM-module-quadric-A-seq} $M$ has a linear resolution
      \begin{equation*}
        \dotsb \to A(-2)^{\oplus s} \to A(-1)^{\oplus s} \to A^{\oplus s} \to M \to 0.
      \end{equation*}
    \item\label{enum:lem-CM-module-quadric-ortho}
      $\Hom_{\gr{A}}(M,A)=0$.
    \item\label{enum:lem-CM-module-quadric-corr}
      Assume further that $A$ is a graded isolated singularity.
      Then the functor $G$ in \ref{lem:CA} induces a bijection between
      the isomorphism classes of indecomposable non-projective modules $M$ in $\CM*{A}$
      with $\dim_{\kk}(M_0)=s$
      and those of simple $C(A)$-modules $N$ with $\dim_{\kk}(N)=s$.
  \end{enumerate}
\end{lemma}

\begin{proof}
  Parts (1) and (2) are proved in
  \cite[Thm.\@ 4.4(1) and Prop.\@ 6.2]{moriNoncommutativeMatrixFactorizations2021}
  and \cite[Prop.\@ 7.8(1)]{moriNoncommutativeMatrixFactorizations2021},
  respectively, using the language of matrix factorizations.
  We provide direct proofs below.

  First, the argument in \cite[Prop.\@ 4.6(1)]{smithNoncommutativeQuadricSurfaces2013}
  (which also applies when $f$ is a normal regular element) shows
  that $M$ admits a minimal free resolution as a graded $S$-module
  \begin{equation}\label{eq:M-S-res-temp}\postdisplaypenalty=10000\relax
    0\to F_1 \xrightarrow{d_1} S^{\oplus s} \to M \to 0
  \end{equation}
  with $\rank{F_1}=s$.

  (1)
  There is an exact sequence of graded $(S, S)$-bimodules
  \[
    0\to S^{\sigma}(-2)
    \xrightarrow{\cdot f}S\to A\to 0.
  \]
  Since $Mf=0$, it follows from the exact sequence
  \[
    0\to\Tor_1^S(M,A)\to M^{\sigma}(-2)
    \xrightarrow{\cdot f}M\to M\to 0
  \]
  that
  $\Tor_1^S(M,A)\cong M^{\sigma}(-2)$.
  Applying $\blank\otimes_S A$ to \ref{eq:M-S-res-temp}
  yields the exact sequence
  \begin{equation*}\label{eq:M-A-pre}
    0\to M^{\sigma}(-2)\to P_1
    \xrightarrow{\partial_1}P_0
    \to M\to 0,
  \end{equation*}
  where $P_0\defeq A^{\oplus s}$ and $P_1\defeq F_1\otimes_S A$ are finite graded free $A$-modules,
  and $\partial_1\defeq d_1\otimes_S \id_A$.
  Then one has $\Omega^2(M)\cong M^{\sigma}(-2)$
  since $M$ has no graded free direct summand.

  (2)
  Since $\ker(\partial_1)=M^{\sigma}(-2)$ has no graded free direct summand,
  $M_A$ has a minimal free resolution
  \begin{equation}\label{eq:M-A-res}
    \dotsb \to P_3\to P_2\xrightarrow{\partial_2} P_1\xrightarrow{\partial_1} P_0\to M\to 0,
  \end{equation}
  where $P_{i+2}\cong P_i^{\sigma}(-2)$ for every $i\geq 0$.
  It suffices to show that $P_1\cong A(-1)^{\oplus s}$.

  Write $P_1=\bigoplus_{i=1}^{s}A(-a_i)$,
  where $a_i\geq 1$ for every $i=1,\dotsc,s$.
  Since $P_2\cong A(-2)^{\oplus s}$,
  we claim that
  \begin{equation}\label{eq:min-res-ineq}
    \max\{a_i\mid 1\leq i\leq s\}<2.
  \end{equation}
  This is a special case of
  \cite[Eq.\@ 3.4]{stephensonGrowthGradedNoetherian1997},
  and we repeat their proof for the convenience of the reader.
  Put $(\blank)^\vee=\gHom_A(\blank,A)$.
  If \ref{eq:min-res-ineq} fails, there is some
  $i_0\in\{1,\dotsc,s\}$ such that $a_{i_0}\geq 2$.
  Since \ref{eq:M-A-res} is minimal, $\partial_2$ sends the generator
  of each summand $A(-2)$ into $\bigoplus_{i\neq i_0}A(-a_i)$, whence
  $\im(\partial_2)\subseteq\bigoplus_{i\neq i_0}A(-a_i)$.
  Hence $\partial_2^\vee(A(a_{i_0}))=0$.
  Since $\gExt_A^1(M,A)=0$, $A(a_{i_0})$
  is in the image of $\partial_1^\vee$.
  But the image of $\partial_1^\vee$ is in
  $\bigoplus_i A_{\geq 1}(a_i)$, a contradiction.
  Therefore \ref{eq:min-res-ineq} holds.
  It follows that $a_i=1$ for every $i$, so $P_1\cong A(-1)^{\oplus s}$.


  (3) By (1) and (2), there is a short exact sequence
  \begin{equation*}
    0 \to M(-2) \to A(-1)^{\oplus s} \to \Omega(M)^{\sigma^{-1}} \to 0.
  \end{equation*}
  Applying $\Hom_{\gr{A}}(\blank, A(-2))$ yields the exact sequence
  \begin{equation*}
    \begin{tikzcd}[cramped,column sep=0.8em]
      \Hom_{\gr{A}}(A(-1)^{\oplus s}, A(-2)) \rar
      & \Hom_{\gr{A}}(M(-2), A(-2)) \rar
      & \Ext_{\gr{A}}^1(\Omega(M)^{\sigma^{-1}}, A(-2)).
    \end{tikzcd}
  \end{equation*}
  The left term is isomorphic to $A_{-1}^{\oplus s}$ and is thus zero,
  while the right term is also zero since $\Omega(M)^{\sigma^{-1}}$ is Cohen--Macaulay.
  It follows that $\Hom_{\gr{A}}(M,A)\cong\Hom_{\gr{A}}(M(-2),A(-2))=0$.

  (4) It follows along the lines of the proof of
  \cite[Prop.\@ 5.4]{smithNoncommutativeQuadricSurfaces2013}.
\end{proof}

\begin{theorem}\label{thm:tilting-hypersurface}
  Let $A$ be a noncommutative isolated quadric hypersurface singularity
  of dimension $d\geq 2$.
  Let $\{M_1,\dotsc,M_n\}$ be a complete set of non-isomorphic indecomposable non-projective modules in $\CM*{A}$.
  Then
  \begin{equation}\label{eq:exceptional-seq}
    \bigl(\pi(A)(-d+2),\pi(A)(-d+3),\dotsc,\pi(A),\pi(M_1),\dotsc,
    \pi(M_n)\bigr)
  \end{equation}
  is a full strong exceptional sequence in $\Dbqgr{A}$.
  Consequently, there is an equivalence of triangulated categories
  \begin{equation}\label{eq:derived-equiv}
    \Dbqgr{A}\simeq \Dbmod{\Lambda},
  \end{equation}
  where $\Lambda \defeq
\End_{\gr{A}}\bigl(A(-d+2)\oplus\dotsb\oplus A(-1)\oplus A
  \oplus M_1\oplus\dotsb\oplus M_n\bigr)$ is the endomorphism algebra, and $\Lambda$
  is finite-dimensional and basic.
\end{theorem}

\begin{proof}
  The derived category $\Dbqgr{A}$ has the Orlov decomposition
  \cite[Thm.\@ 16]{orlovDerivedCategoriesCoherent2009}
  \begin{equation}\label{eq:orlov-decomp}
    \Dbqgr{A}=\langle\pi(A)(-d+2),\pi(A)(-d+3),\dotsc,\pi(A),
    \Phi(\sCM{A})\rangle,
  \end{equation}
  where $\Phi\colon \sCM{A}\to \Dbqgr{A}$ is a fully faithful functor.
  From the construction of $\Phi$
  (see also \cite[Lem.\@ 2.6]{ueyamaDerivedCategoriesSkew2023}),
  we have $\Phi(M_i)\cong\pi(M_i)$ for $1\leq i\leq n$
  by \ref{enum:lem-CM-module-quadric-ortho}.
  By \ref{lem:CA,thm:quadric-iso-sing},
  $C(A)$ is semisimple,
  and the triangulated equivalence
  $F\colon \sCM{A}\simeq \Dbmod{C(A)}$
  sends $M_1,\dotsc, M_n$ to stalk complexes of the simple $C(A)$-modules
  concentrated in degree 0.
  Therefore, $(M_1, \dotsc, M_n)$
  is a full strong exceptional sequence
  in $\sCM{A}$.
  It follows from \ref{eq:orlov-decomp}
  that \ref{eq:exceptional-seq} is a full exceptional sequence.
  To conclude that it is strong,
  it remains to show
  \begin{align}
    \label{eq:tilting-hypersurface-1}
    \Ext_{\qgr{A}}^q\bigl(\pi(A(-j)),\pi(M_i)\bigr)=0,
    &\quad \text{for }1\leq i\leq n, 0\leq j\leq d-2,q>0,
    \text{ and}\\
    \label{eq:tilting-hypersurface-2}
    \Ext_{\qgr{A}}^q\bigl(\pi(A(-j)),\pi(A(-i))\bigr)=0,
    &\quad \text{for }0\leq i\leq j\leq d-2,q>0.
  \end{align}
  Since $A(-j)$ is projective,
  the left-hand side of \ref{eq:tilting-hypersurface-1} is zero
  for $q\neq d-1$
  by \ref{lem:Serre-duality}(1) and (3).
  On the other hand, since $M_i$ is generated in degree 0 and $-j-d+1<0$,
  we obtain
  \begin{align*}
    \Ext_{\qgr{A}}^{d-1}\bigl(\pi(A(-j)),\pi(M_i)\bigr)&\cong
    D\bigl(\Hom_{\gr{A}}(M_i,A^{\mu^{-1}}(-j-d+1))\bigr)=0.
  \end{align*}
  Therefore, \ref{eq:tilting-hypersurface-1} holds.
  The same argument establishes \ref{eq:tilting-hypersurface-2},
  since $-j-d+1+i<0$ for $0\leq i\leq j\leq d-2$.

  Let $T=A(-d+2)\oplus\dotsb\oplus A\oplus M_1\oplus\dotsb\oplus M_n$.
  Since $\depth{T}=d\geq 2$, $\Lambda=\End_{\gr{A}}(T)\cong\End_{\qgr{A}}(\pi(T))$.
  By \cite[Cor.\@ 4.3]{yekutieliSerreDualityNoncommutative1997} and
  \cite[Lem.\@ 2.5]{moriMcKaytypeCorrespondenceASregular2013},
  $\qgr{A}$ is Hom-finite, and thus $\Lambda$ is finite-dimensional.
  Since $\QGr{A}$ has enough injective objects,
  the argument in
  \cite[Thm.\@ 6.2]{bondalRepresentationAssociativeAlgebras1990},
  together with the subsequent remark, shows that $\RHom_{\QGr{A}}(\pi(T),\blank)$
  induces the equivalence \ref{eq:derived-equiv}.
\end{proof}


\section{Cluster tilting for
noncommutative isolated quadric hypersurface singularities}\label{sec:ct-hyper}

In this section, we investigate the existence of graded $(d-1)$-cluster tilting modules
over noncommutative isolated quadric hypersurface singularities
of dimension $d\geq 2$.

The notion of an $n$-cluster tilting module was first introduced under the name of an $(n-1)$-maximal orthogonal module in
\cite[Def.\@ 2.4]{iyamaAuslanderCorrespondence2007}.
It is central to the study of
noncommutative (crepant) resolutions of singularities
\cite[Thm.\@ 5.2.1]{iyamaAuslanderCorrespondence2007}.

\begin{definition}
  Let $A$ be a noetherian AS-Gorenstein algebra
  and let $n\geq 1$.
  A \dfn{graded $n$-cluster tilting module} is a module $M\in\CM{A}$ such that
  \begin{align*}
    \add(M)&=\{X\in\CM{A}\mid
    \gExt_A^i(M,X)=0\text{ for }0<i<n\}\\
    &=\{X\in\CM{A}\mid
    \gExt_A^i(X,M)=0\text{ for }0<i<n\}.
  \end{align*}
\end{definition}

We recall the definition of a noncommutative resolution of an AS-Gorenstein isolated singularity given in \cite[Def.\@ 6.1]{liNoncommutativeResolutionsASGorenstein2026}.

\begin{definition}
  \label{dfn:NCCR}
  Let $A$ be an AS-Gorenstein isolated singularity
  of dimension $d\geq 2$.
  A \dfn{noncommutative resolution} of $A$
  is a noetherian $\ZZ$-graded algebra $B$ of global dimension $d$,
  of the form $B=\gEnd_A(M)$,
  where $M_A$ is a Cohen--Macaulay generator
  and $B_B$ is Cohen--Macaulay.
\end{definition}

\begin{theorem}
  [{{\cite[Thm.\@ 6.6]{liNoncommutativeResolutionsASGorenstein2026}}}]
  \label{thm:NCCR-ct}
  Let $A$ be an AS-Gorenstein
  isolated singularity of dimension $d\geq 2$,
  and let $M$ be a Cohen--Macaulay generator of $A$.
  Then $\gEnd_A(M)$ is a noncommutative resolution of $A$
  if and only if $M$ is a graded $(d-1)$-cluster tilting module.
\end{theorem}


We first consider the graded endomorphism algebra
of an additive generator of $\CM{A}$,
as studied in
\cite{ueyamaNoncommutativeGradedAlgebras2015,
chanMcKayCorrespondenceSemisimple2019}.
These results extend Leuschke's result on commutative local rings
\cite{leuschkeEndomorphismRingsFinite2007}.

\begin{proposition}\label{prop:high-dim-good}
  Suppose that $A$ is an AS-Gorenstein isolated singularity
  of dimension $d\geq 2$
  with an additive generator $M$ of $\CM{A}$.
  Let $B\defeq\gEnd_A(M)$. Then
  \begin{enumerate}[ref={\ref{prop:high-dim-good}(\arabic*)}]
    \item $B$ is noetherian;
    \item\label{enum:prop-high-dim-good-qgr}
      there is an equivalence
      $(\qgr A, (1))\simeq (\qgr B, (1))$
      sending $\pi(M)$ to $\pi(B)$; and
    \item\label{enum:prop-high-dim-good-gldim} $\gldim(B)=d$.
  \end{enumerate}
\end{proposition}

\begin{proof}
  Statements (1) and (2) are shown in
  \cite[Thm.\@ 3.8(2)]{ueyamaNoncommutativeGradedAlgebras2015}.
  Statement (3) is a special case of
  \cite[Thm.\@ 5.4(2)]{chanMcKayCorrespondenceSemisimple2019}.
\end{proof}


In \ref{prop:high-dim-good}, $B$ must be Cohen--Macaulay to qualify as a noncommutative resolution.
\begin{theorem}\label{thm:high-dim-bad}
  Let $A$ be a noetherian
  AS-Gorenstein algebra of dimension $d\geq 3$
  with finite graded Cohen--Macaulay representation type.
  Let $M$ be an additive generator of $\CM{A}$,
  and let $B=\gEnd_A(M)$. If $A$ is not AS-regular, then
  \begin{enumerate}[ref={\ref{thm:high-dim-bad}(\arabic*)}]
    \item $\gExt_A^1(M, M)\neq 0$; and
    \item\label{enum:high-dim-bad-depth} $\depth{B}=2<d$.
  \end{enumerate}
\end{theorem}

\begin{proof}
  (1) Let $N$ denote the $d$-th syzygy of the trivial module $\kk$.
  Then $N\in\CM{A}=\add(M)$.
  Suppose, to the contrary, that $\gExt_A^1(M, M)=0$.
  Then $\gExt_A^{d+1}(\kk, M)\cong\gExt_A^1(N, M)=0$ by dimension shifting.
  Let $I^{d+1}$ denote the $(d+1)$-th term in the minimal injective resolution of $M$.
  Since $\soc{I^{d+1}}=\gExt_A^{d+1}(\kk, M)=0$,
  $I^{d+1}$ is torsion-free.
  On the other hand,
  since $A$ has finite graded Cohen--Macaulay representation type,
  it follows from \cite[Thm.\@ 3.4]{ueyamaNoncommutativeGradedAlgebras2015}
  and \ref{eq:iso-sing} that $\gldim(\QGr{A})=d-1$, and thus $I^{d+1}$ is torsion.
  Consequently, $I^{d+1}=0$, so $\idim(M)\leq d$.
  Thus, every Cohen--Macaulay $A$-module has injective dimension $\leq d$.

  For any $X\in\gr{A}$, by Auslander--Buchweitz approximation
  \cite[Prop.\@ 8.7]{liCommonlyGradedAlgebras2025}
  (see also \cite[Prop.\@ 5.3]{moriHomologicalPropertiesBalanced2003}),   there is a short exact sequence
  $0\to Z\to Y\to X\to 0$,
  where $Y$ is Cohen--Macaulay and $\idim(Z)\leq d$.
  Since $\idim(Y)\leq d$, it follows that $\idim(X)\leq d$.
  Therefore, $\gldim(A)\leq d$, contrary to the assumption that
  $A$ is not AS-regular.

  (2) First, $\depth{B}\geq 2$ by
  \cite[Cor.\@ 3.8]{liNoncommutativeResolutionsASGorenstein2026}.
  For the reverse inequality,
  \ref{enum:prop-high-dim-good-qgr} and
  \cite[Lem.\@ 5.5(2) and Prop.\@ 3.3(1)]{%
  liNoncommutativeResolutionsASGorenstein2026}
  yield the isomorphisms
  \begin{equation*}
0\neq\gExt_A^1(M,M)\cong\gExt_{\QGr{A}}^1(\pi(M),\pi(M))
    \cong\gExt_{\QGr{B}}^1(\pi(B),\pi(B))
    \cong\uR^2\Gamma_B(B),
  \end{equation*}
  which imply that $\depth{B}\leq 2$.
\end{proof}

Thus, for $d\geq 3$, an additive generator of $\CM{A}$
does not give a noncommutative resolution unless $A$ is AS-regular.
We therefore ask whether a noncommutative resolution can be obtained
from another Cohen--Macaulay generator of $A$.

In the commutative local setting, a result similar to \ref{enum:high-dim-bad-depth} was proved by Graham J.\@ Leuschke in a 2026 note on his website, as an erratum to \cite[Thm.\@ 6]{leuschkeEndomorphismRingsFinite2007}. 


For the remainder of this section, let $A=S\div(f)$ be
a noncommutative isolated quadric hypersurface singularity
of dimension $d\geq 2$.
Recall that $A$ has finite graded Cohen--Macaulay representation type
by \ref{thm:quadric-iso-sing}.

By \ref{lem:CM-module-quadric},
$\rho\defeq\Omega(\blank)(1)$ defines an autoequivalence of $\sCM*{A}$ and therefore induces a permutation of the finite set $\ind(\sCM*{A})$.
We regard the objects of $\ind(\sCM*{A})$ as indecomposable non-projective modules in $\CM*{A}$.
The following lemma is used to verify the cluster tilting condition.

\begin{lemma}\label{lem:Ext-indec}
  Let $X, Y\in\ind(\sCM*{A})$ and $q>0$.
  Then $\gExt_A^q(X,Y)\neq 0$ if and only if $Y\cong\rho^q(X)$.
\end{lemma}

\begin{proof}
  Recall that there is a triangulated equivalence
  $F\colon \sCM{A}\simeq \Dbmod{C(A)}$ (\ref{enum:lem-CA-F})
  and that $C(A)$ is semisimple (\ref{thm:quadric-iso-sing}).
  By \ref{enum:lem-CA-G}, $F$ sends the objects of
  $\ind(\sCM*{A})$ to stalk complexes of
  simple $C(A)$-modules concentrated in degree 0.
  Therefore, for $X,Y\in\ind(\sCM*{A})$, $q>0$ and $t\in\ZZ$,
  \[\Ext_{\gr{A}}^q(X,Y(t))
  \cong\Hom_{\sCM{A}}(X, \rho^t(Y)[t+q])
  \cong\Ext^{t+q}_{C(A)}\bigl(F(X), F(\rho^t(Y))\bigr),\]
  which is nonzero if and only if $X\cong \rho^t(Y)$ and $t=-q$.
  Taking the direct sum over all $t$
  shows that $\gExt_A^q(X,Y)\neq 0$ if and only if $Y\cong\rho^q(X)$.
\end{proof}

We now use the permutation $\rho$ to characterize when $A$ admits a graded $(d-1)$-cluster tilting module.

\begin{proposition}\label{prop:ct-decomp}
  Let $A$ be a noncommutative isolated quadric hypersurface singularity
  of dimension $d\geq 2$.
  The following conditions are equivalent.
  \begin{enumerate}
    \item $A$ has a graded $(d-1)$-cluster tilting module.
    \item There exists a subset $I\subseteq\ind(\sCM*{A})$ such that
      \[\ind(\sCM*{A})=I\sqcup \rho(I) \sqcup\dotsb\sqcup \rho^{d-2}(I),\]
      where $\sqcup$ denotes disjoint union and $\rho$ is defined before \ref{lem:Ext-indec}.
    \item The cardinality of every $\rho$-orbit of $\ind(\sCM*{A})$
      is divisible by $d-1$.
  \end{enumerate}
\end{proposition}

\begin{proof}
  For a subset $I$ of $\ind(\sCM*{A})$,
  we use $M_I$ to denote $A\oplus\bigoplus_{M\in I}M$.

  \textit{(1)$\Rightarrow$(2).}
  Let $M$ be a graded $(d-1)$-cluster tilting module.
  Let $I=\add(M)\cap\ind(\sCM*{A})$.
  Since $A\in\add(M)$, it follows from \ref{lem:CM-shift}
  that $\add(M)=\add(M_I)$.
  Hence, $M_I$ is also a graded $(d-1)$-cluster tilting module.

  Since $M_I$ is $(d-1)$-rigid, that is, $\gExt_A^q(M_I, M_I)=0$ for $0<q<d-1$,
  it follows from \ref{lem:Ext-indec} that $I\cap \rho^q(I)=\emptyset$.
  Consequently, $\rho^p(I)\cap \rho^q(I)=\emptyset$ for $0<p<q<d-1$.

  Suppose $Y\in \ind(\sCM*{A})\setminus I$.
  Since $M_I$ is graded $(d-1)$-cluster tilting,
  there exists $q$ with $0<q<d-1$ such that $\gExt_A^q(M_I, Y)\neq 0$.
  Hence, $Y\in\rho^q(I)$ by \ref{lem:Ext-indec}.

  \textit{(2)$\Rightarrow$(1).}
  Let $X\in\ind(\sCM*{A})$.
  By \ref{lem:Ext-indec}, $\gExt_A^i(M_I, X)=0$ for $0<i<d-1$ if and only if $X\in I$.
  Moreover, $\gExt_A^i(M_I,A)=0$ for $i>0$ since $M_I$ is Cohen--Macaulay.
  It follows from \ref{lem:CM-shift} that
  \[\add(M_I)=\{X\in\CM{A}\mid
  \gExt_A^i(M_I,X)=0\text{ for }0< i< d-1\},\]
  and thus $M_I$ is graded $(d-1)$-cluster tilting
  by \cite[Thm.\@ 6.6]{liNoncommutativeResolutionsASGorenstein2026}.

  \textit{(2)$\Rightarrow$(3).}
  Let $O$ be a $\rho$-orbit. Since $\rho(O)=O$, we have
  \[
    O=(O\cap I)\sqcup \rho(O\cap I) \sqcup\dotsb\sqcup \rho^{d-2}(O\cap I).
  \]
  Hence, $\abs{O}$ is divisible by $d-1$.

  \textit{(3)$\Rightarrow$(2).}
  For each $\rho$-orbit $O$ of cardinality $n(d-1)$, fix an element $X\in O$
  and let $I_O$ denote the set $\{\rho^{i(d-1)}(X)\mid 0\leq i<n\}$,
  which has cardinality $n$.
  Then $O=I_O\sqcup \rho(I_O)\sqcup \dotsb \sqcup \rho^{d-2}(I_O)$.
  Let $I$ be the union of the sets $I_O$ over all $\rho$-orbits $O$.
  It follows that
  \[\ind(\sCM*{A})=I\sqcup \rho(I) \sqcup\dotsb\sqcup \rho^{d-2}(I).\qedhere\]
\end{proof}

\begin{corollary}\label{cor:central-no-NCCR}
  If $A$ is a noncommutative isolated central quadric hypersurface singularity
  of dimension $d\geq 4$,
  then $A$ has no graded $(d-1)$-cluster tilting module.
\end{corollary}

\begin{proof}
  Since $f$ is central, $\rho^2$ acts as the identity on
  $\ind(\sCM*{A})$ by \ref{lem:CM-module-quadric}.
  Therefore, every $\rho$-orbit has cardinality 1 or 2,
  and condition (3) of \ref{prop:ct-decomp} cannot hold when $d\geq 4$.
\end{proof}

\begin{remark}\label{rem:comm-no-NCCR}
  Let $S$ be a quantum polynomial algebra of dimension $d+1$, and let $f\in S_2$ be a normal regular element.
  Suppose that $A=S/(f)$ is commutative.
  Since $A$ is generated in degree 1
  and has Hilbert series $(1-t)^{-d}(1+t)$,
  there is a graded isomorphism $\kk[x_0,\dotsc,x_d]/(g)\to A$
  for some nonzero quadratic polynomial $g$.
  It follows from \ref{cor:central-no-NCCR}
  that if $A$ is a graded isolated singularity and $d\geq 4$,
  then $A$ does not admit a graded $(d-1)$-cluster tilting module.
\end{remark}

\begin{remark}\label{rem:standard}
  Suppose that $A=S/(f)$ is
  a noncommutative isolated quadric hypersurface singularity of dimension 3
  and that $\ind(\sCM*{A})$ consists of exactly two elements, $X$ and $Y$.
  By \ref{prop:ct-decomp},
  $A$ has a graded 2-cluster tilting module
  if and only if $\rho(X)\cong Y$ and $\rho(Y)\cong X$. If this condition holds and $f$ is irreducible in $S$, then $A=S/(f)$ is called \dfn{standard} in \cite[Def.\@ 5.4]{moriCategoricalCharacterizationQuantum2021} and
  \cite[Def.\@ 3.7]{moriIrreducibleNoncommutativeQuadrics2026}.
\end{remark}

Combining \ref{rem:standard} with \cite[Thm.\@ 5.17(2)]{moriIrreducibleNoncommutativeQuadrics2026},
we obtain the following corollary.
For unexplained terminology, see \cite{moriIrreducibleNoncommutativeQuadrics2026}.

\begin{corollary}
  Assume $\charOp(\kk)=0$.
  Let $S$ be a 4-dimensional non-degenerate Sklyanin algebra,
  and let $0\neq f\in S_2$ be central.
  If $A=S/(f)$ is a graded isolated singularity,
  then $A$ admits a graded 2-cluster tilting module.
\end{corollary}

We now characterize the values of $d$ for which
there is a noncommutative isolated quadric hypersurface singularity of dimension $d$
with a graded $(d-1)$-cluster tilting module.


\begin{theorem}\label{thm:ct-dimension}
  The following statements are equivalent.
 \begin{enumerate}
     \item There exists a noncommutative isolated quadric hypersurface singularity
  of dimension $d$ with a graded $(d-1)$-cluster tilting module, where $d\geq 2$.
  \item $d=2^n+1$ for some integer $n\geq 0$.
  \end{enumerate}
\end{theorem}

\begin{proof} \textit{(1)$\Rightarrow$(2).}
  Suppose $A$ is a noncommutative isolated quadric hypersurface singularity
  of dimension $d$ that admits a graded $(d-1)$-cluster tilting module.
  Then $C(A)$ is semisimple (\ref{thm:quadric-iso-sing}),
  and $\dim_{\kk}(C(A))=2^d$ (\ref{enum:lem-CA-dim}).
  By \ref{enum:lem-CM-module-quadric-A-seq},
  $\dim_{\kk}(M_0)=\dim_{\kk}(\rho(M)_0)$ for any $M\in\ind(\sCM*{A})$.
  Fix $I$ as in \ref{prop:ct-decomp}.
  Let $\simp{C(A)}$ be the set of isomorphism classes
  of simple $C(A)$-modules.
  Since $C(A)$ is semisimple and $\kk$ is algebraically closed, it follows from
  \ref{enum:lem-CM-module-quadric-corr} that
  \begin{align*}
    2^d = \dim_{\kk}(C(A)) &= \sum_{V\in\simp{C(A)}} \dim_{\kk}(V)^2\\
    &=\sum_{M\in\ind(\sCM*{A})} \dim_{\kk}(M_0)^2\\
    &=(d-1)\sum_{M\in I} \dim_{\kk}(M_0)^2.
  \end{align*}
  Thus $d-1$ divides $2^d$, so $d=2^n+1$ for some $n\geq 0$.

 \textit{(2)$\Rightarrow$(1).} Suppose that $d=2^n+1$ for some $n\geq 0$.
  If $n=0$, then $d=2$.
  A graded 1-cluster tilting module is precisely
  an additive generator of $\CM{A}$; hence, such a module exists if and only if
  $A$ has finite graded Cohen--Macaulay representation type.
  Thus, every noncommutative isolated quadric surface singularity
  gives an example by \ref{thm:quadric-iso-sing}.
  For $n\geq 1$, a noncommutative isolated quadric hypersurface singularity
  of dimension $2^n+1$ with a graded $(d-1)$-cluster tilting module
  will be constructed in \ref{ex:ct}, completing the proof.
\end{proof}

In the setting of \ref{rem:standard},
Mori, Ueyama, and Wu studied standardness under Zhang twists 
\cite[Sec.\@ 4]{moriIrreducibleNoncommutativeQuadrics2026}.
Inspired by their approach, we construct the following example
as a graded twist of a skew quadric hypersurface studied in \cite{ueyamaDerivedCategoriesSkew2023}.

\begin{example}\label{ex:ct}
  Suppose $n\geq 1$ and $d=2^n+1$.
  Let
  \[S=\kk\langle x_1,\dotsc, x_{d+1}\rangle
  \div (x_ix_j+x_jx_i)_{1\leq i<j\leq d+1},\]
  let $f=x_1^2+x_2^2+\dotsb+x_{d+1}^2$, which is central and regular,
  and let $A=S\div(f)$.
  By \cite[Thm.\@ 1.4]{higashitaniCombinatorialStudyStable2022}, $A$ is a graded isolated singularity.
  Moreover, $\Omega(M)(1)\cong M$ for every $M\in\ind(\sCM*{A})$
  by \cite[Ex.\@ 3.26]{ueyamaDerivedCategoriesSkew2023}.

  Let $\Tau$ be the graded automorphism of $S$ defined by
  \[
    \Tau(x_i)=
    \begin{cases}
      x_{i+1}, &i<2^{n-1},\\
      -x_1, &i=2^{n-1},\\
      x_i, &i>2^{n-1}.
    \end{cases}
  \]
  Then $\Tau^{2^{n-1}}$ sends $x_i$ to $-x_i$ for $i\leq 2^{n-1}$,
  and fixes $x_i$ for $i>2^{n-1}$.
  Consequently, $\Tau$ has order $2^n$.
  Since $\Tau(f)=f$, $\Tau$ induces a graded automorphism $\tau$ of $A$.

  Define $A'=A_\tau$.
  By \ref{enum:prop-zhang-twist-noeth}, \ref{enum:prop-zhang-twist-Koszul}, and \ref{lem:twist-normal-regular}, $A'\cong S_\Tau\div(f_\Tau)$ is also
  a noncommutative isolated quadric hypersurface singularity
  of dimension $d$.
  We show that $A'$ admits a graded $(d-1)$-cluster tilting module.

  Let $\rho'\defeq\Omega_{A'}(\blank)(1)$.
  By \ref{enum:prop-zhang-twist-CM},
  the equivalence $F\colon \gr{A}\to\gr{A_\tau}$ defined by $M\mapsto M_\tau$
  preserves Cohen--Macaulay modules and induces a bijection
  $\ind(\sCM*{A})\cong\ind(\sCM*{A_\tau})$.
  There is also an isomorphism $F(M)(1)\cong F(M^\tau(1))$ induced by the identity map on the underlying vector space.
  It follows that for $M\in\ind(\sCM*{A})$,
  \[\Omega(F(M))(1)\cong F(\Omega(M))(1)\cong F(\Omega(M^\tau)(1))
  \cong F(M^\tau),\]
  and hence $\rho'(F(M))=\Omega(F(M))(1)\cong F(M^\tau)$.

  By \ref{prop:ct-decomp}, it suffices to show that every $\rho'$-orbit
  on $\ind(\sCM*{A'})$ has cardinality $2^n=d-1$.
  Fix $M\in\ind(\sCM*{A})$.
  By \cite[Ex.\@ 3.26]{ueyamaDerivedCategoriesSkew2023},
  $M\cong A/gA$ for some $g=x_1+\sum_{i=2}^{d+1}a_ix_i$
  with $a_i\in\{\pm 1\}$.


  Let $\Theta=\Tau^{2^{n-1}}$ and $\theta=\tau^{2^{n-1}}$.
  We claim that $M\ncong M^\theta$.
  Suppose that $\phi\colon M\to M^\theta$ is a nonzero graded map.
  Let $e\in M_0$ and $e'\in M^\theta_0$ be nonzero.
  Then there exists $\lambda\neq 0$ such that $\phi(e)=\lambda e'$.
  Since $eg=0$,
  we have $0=\phi(eg)=\phi(e)\cdot^\theta g=\lambda e'\theta(g)$.
  Thus $\theta(g)\in \kk g\subseteq A_1$, and hence $\Theta(g)\in \kk g\subseteq S_1$.
  However, $\Theta(x_1)=-x_1$ and $\Theta(x_d)=x_d$,
  while the coefficients of $x_1$ and $x_d$ in $g$ are nonzero.
  This is impossible since $\charOp(\kk)\neq 2$.

  Consequently, $\rho'^{2^{n-1}}(F(M))\cong F(M^\theta)\ncong F(M)$.
  Since $\tau^{2^n}=\id_A$, every $\rho'$-orbit on $\ind(\sCM*{A'})$ has cardinality dividing $2^n$.
  Any $\rho'$-orbit of cardinality less than $2^n$ would have
  cardinality dividing $2^{n-1}$, contradicting the preceding claim.
  Therefore, every $\rho'$-orbit on $\ind(\sCM*{A'})$ has cardinality $2^n=d-1$.
  By \ref{prop:ct-decomp}, $A'$ admits a graded $(d-1)$-cluster tilting module.
\end{example}


\section{Cluster tilting for
noncommutative isolated quadric surface singularities}\label{sec:ct-surface}

In this section, we construct noncommutative resolutions of noncommutative isolated quadric surface singularities
by using graded cluster tilting modules.
These resolutions are twisted preprojective algebras of certain types.
We establish a correspondence between these twisted preprojective algebras and noncommutative isolated quadric surface singularities (\ref{thm:classif-corresp}).
By classifying these twisted preprojective algebras, we obtain a complete classification of the singularities up to graded isomorphism (\ref{thm:A-classif}).

\subsection{Graded preprojective algebras and their twists}

We recall the construction of graded preprojective algebras
and collect their homological and ring-theoretic properties.
We also introduce their graded twists
and characterize twisted preprojective algebras in terms of the graded skew CY property.
For background on quivers, we refer to
\cite[Chap.\@ II]{assemElementsRepresentationTheory2006}.

A \dfn{quiver} $Q$ consists of
the set $Q_0$ of vertices of $Q$,
the set $Q_1$ of arrows of $Q$,
and two maps $s,t\colon Q_1\to Q_0$
sending each arrow to its source and target, respectively.

Let $Q$ be a quiver. Then $Q$ is \dfn{finite} if both $Q_0$ and $Q_1$ are finite sets; $Q$ is \dfn{connected} if the underlying undirected diagram is connected; and $Q$ is \dfn{acyclic} if $Q$ has no oriented cycles.

Let $Q$ be a finite quiver, and let $\kk Q$ be its path algebra.
We identify each vertex with its corresponding trivial path.
We adopt the convention that arrows in $\kk Q$ are composed from left to right.

A \dfn{weight} on a quiver $Q$ is a map $w\colon Q_1\to\NN$.
A quiver equipped with a weight is called a \dfn{weighted quiver}.
If $Q$ is a weighted quiver with weight $w$,
then $\kk Q$ is an $\NN$-graded algebra
in which each $e\in Q_0$ has degree 0
and each $\alpha\in Q_1$ has degree $w(\alpha)$.

\begin{definition}\label{dfn:preproj-alg}
  Let $Q$ be a finite quiver.
  The \dfn{double quiver} $Q^\sharp$ of $Q$
  is obtained by keeping the vertex set
  and adding a new arrow $\alpha^*$ from $t(\alpha)$ to $s(\alpha)$
  for each arrow $\alpha\in Q_1$.
  The \dfn{mesh relation} $\rho_e$ at vertex $e$
  is the element of $\kk Q^\sharp$ given by
  \begin{equation*}
    \rho_e\defeq\sum_{\substack{\alpha\in Q_1\\s(\alpha)=e}}\alpha\alpha^*
    -\sum_{\substack{\alpha\in Q_1\\t(\alpha)=e}}\alpha^*\alpha.
  \end{equation*}
  The total sum $\rho\defeq\sum_{e\in Q_0}\rho_e$
  is called the \dfn{mesh relation} of $Q^\sharp$.
  The \dfn{preprojective algebra} $\Pi(\kk Q)$
  is defined as the quiver algebra
  \begin{equation*}
    \Pi(\kk Q)\defeq \kk Q^\sharp\div(\rho)
    = \kk Q^\sharp\div
    \bigl(\rho_e\mid e\in Q_0\bigr).
  \end{equation*}
\end{definition}

The double quiver $Q^\sharp$ has a natural weight $w$ defined by $w(\alpha)=0$ and $w(\alpha^*)=1$ for every $\alpha\in Q_1$.
This weight induces a grading on $\Pi(\kk Q)$ such that $\Pi(\kk Q)_0=\kk Q$.
With this grading,
$\Pi(\kk Q)$ is referred to as the \dfn{graded preprojective algebra}
in this paper.


A connected finite acyclic quiver $Q$ is called a \dfn{Dynkin quiver} or a \dfn{Euclidean quiver}
if its underlying unoriented graph is a Dynkin or Euclidean graph, respectively
(see \cite[Sec.\@ VII.2]{assemElementsRepresentationTheory2006}).
If $Q$ is neither Dynkin nor Euclidean, it is called a \dfn{wild quiver}.
The following proposition summarizes properties of graded preprojective algebras with references.

\begin{proposition}\label{prop:preproj}
  Let $Q$ be a connected finite acyclic non-Dynkin quiver, and
  let $\Pi=\Pi(\kk Q)$ denote the graded preprojective algebra.
  \begin{enumerate}[ref={\ref{prop:preproj}(\arabic*)}]
    \item $\Pi$ is prime and graded coherent.
    \item\label{enum:prop-preproj-res}
      Every graded simple $\Pi$-module $S_e$ (concentrated in degree 0)
      at a vertex $e\in Q_0$
      has the minimal projective resolution
      \begin{equation*}
        0\to e\Pi(-1) \xrightarrow{p_1}
        \bigoplus_{\substack{\alpha\in Q_1\\s(\alpha)=e}}t(\alpha)\Pi
        \oplus
        \bigoplus_{\substack{\alpha\in Q_1\\s(\alpha^*)=e}}t(\alpha^*)\Pi(-1)
        \xrightarrow{p_0} e\Pi \to S_e \to 0.
      \end{equation*}
    \item\label{enum:prop-preproj-CY}
      $\Pi$ is graded CY of dimension 2
      with Gorenstein parameter 1.
    \item $\Pi$ is generalized Koszul in the sense of
      \cite[Def.\@ 3.6]{liGeneralizedKoszulAlgebra2023}.
  \end{enumerate}
  Moreover, if $Q$ is Euclidean,
  then the following hold.
  \begin{enumerate}[resume*]
    \item\label{enum:prop-preproj-noeth}
      $\Pi$ is a finite module over its noetherian centre.
      In particular, $\Pi$ is noetherian and satisfies a polynomial identity.
    \item $\GKdim(\Pi)=\Kdim(\Pi)=2$.
  \end{enumerate}
  If $Q$ is wild, then the following hold.
  \begin{enumerate}[label={(\arabic*')},start=5]
    \item $\Pi$ is not noetherian on either side.
    \item $\Pi$ has exponential growth and, in particular, $\GKdim(\Pi)=\infty$.
  \end{enumerate}
\end{proposition}

\begin{proof}
  (1) $\Pi$ is prime by \cite[Prop.\@ 3.4]{baerHomologicalPropertiesWild1986}
    and \cite[Lem.\@ 6.1]{baerPreprojectiveAlgebraTame1987}.
    The graded coherence of $\Pi$ is established in
    \cite[Prop.\@ 5.1 and Thm.\@ 3.7(2)]{minamotoAmplenessTwosidedTilting2012}.


  (2) The resolution was constructed in the proof of \cite[Prop.\@ 7.1(2)]{reyesGrowthGradedTwisted2019}.

  (3) A proof using Keller's derived CY completions
    \cite[Thm.\@ 4.8]{kellerDeformedCalabiYau2011}
    can be found in
    \cite[Thm.\@ 4.36]{herschend$n$representationInfiniteAlgebras2014}.
    See also
    \cite[Lem.\@ 13]{martinez-villaSkewGroupAlgebras2001}
    and
    \cite[Thm.\@ 4.2]{minamotoStructureASGorensteinAlgebras2011}.

  (4) The associated graded algebra of $\Pi$
    with respect to its graded Jacobson radical
    is isomorphic to the standard graded preprojective algebra,
    which is Koszul by a regrading of (2) or by
    \cite[Thm.\@ 1.9]{martinez-villaApplicationsKoszulAlgebras1996}.
    Thus, $\Pi$ is generalized Koszul by
    \cite[Thm.\@ 5.3(2)]{liGeneralizedKoszulAlgebra2023}.

  (5) It is proved in \cite[Thm.\@ C]{braunFinitelyGeneratedPI1994}.


  (5') By \cite[Thm.\@ 4.6]{lenzingHomologicalTransferFinitely1983},
    the category $\cP$ of preprojective $\kk Q$-modules is not left noetherian
    in the sense of \cite[p.\@ 19]{mitchellRingsSeveralObjects1972}.
    Consequently, the category of covariant additive functors $\cP\to \Ab$,
    which is equivalent to the category of $\NN$-graded left modules over $\Pi$ by
    \cite[Prop.\@ 3.6]{baerPreprojectiveAlgebraTame1987},
    is not locally noetherian.
    It follows that $\Pi$ is not left noetherian; neither is $\Pi^\op$,
    since $Q^\op$ also has wild representation type.

  (6) and (6') They are shown in
    \cite[Cor.\@ and Thm.]{dlabEigenvaluesCoxeterTransformations1981}.
\end{proof}

We next consider twists of graded preprojective algebras.

\begin{definition}\label{dfn:twisted-preproj}
  Let $Q$ be a finite acyclic non-Dynkin quiver,
  and let $\Pi$ be its graded preprojective algebra.
  If $\tau$ is a graded automorphism of $\Pi$ that permutes the vertices of $Q$,
  then the twisted algebra $\Pi_{\tau}$
  is called the \dfn{$\tau$-twisted preprojective algebra of type $Q$}.
\end{definition}

\begin{remark}\label{rem:graded-iso-idem}
  Let $\phi\colon B\to B'$ be a graded isomorphism of $\NN$-graded algebras. By assumption, $B_0$ and $B'_0$ are finite-dimensional.
  Let $I\subseteq B_0$ and $I'\subseteq B'_0$ be sets
  of complete orthogonal primitive idempotents of $B_0$ and $B'_0$, respectively.
  Since $\phi(I)$ and $I'$ are sets of complete orthogonal primitive idempotents in $B'_0$,
  there exists a unit $u\in B'_0$ such that $u\phi(I)u^{-1}=I'$
  by \cite[Exer.\@ 21.17]{lamFirstCourseNoncommutative2001}.
  Thus, $b\mapsto u\phi(b)u^{-1}$ is a graded isomorphism from $B$ to $B'$
  sending $I$ to $I'$.

  In particular, given a graded isomorphism $\phi\colon B\to B'$
  between two twisted preprojective algebras of the same type $Q$,
  there exists an inner automorphism $\eta$ of $B'$
  such that $\eta\phi$ permutes the vertices.
\end{remark}

The following proposition, based on \cite{minamotoStructureASGorensteinAlgebras2011}, characterizes twisted preprojective algebras in terms of the graded skew CY property.

\begin{proposition}\label{cor:twisted-preproj}
  Let $B$ be an $\NN$-graded algebra such that $B_0$ is basic and indecomposable as an algebra. 
  The following statements are equivalent.
  \begin{enumerate}[ref={\ref{cor:twisted-preproj}(\arabic*)}]
    \item $B$ is graded skew CY of dimension 2 with Gorenstein parameter 1.
    \item $B$ is graded isomorphic to a twisted preprojective algebra of type $Q$
      for some connected finite acyclic non-Dynkin quiver $Q$.
  \end{enumerate}
  Moreover, if these conditions hold, $B_0\cong\kk Q$,
  and $B$ is noetherian if and only if $Q$ is Euclidean.
\end{proposition}

\begin{proof}
  \textit{(1)$\Rightarrow$(2).}
  Let $\Lambda=B_0$.
  It follows from \ref{thm:CY-equiv} and \cite[Thm.\@ 4.12(1)]{minamotoStructureASGorensteinAlgebras2011}
  that $\Lambda$ is hereditary of infinite representation type.
  Therefore, we may identify $\Lambda$ with the path algebra $\kk Q$
  for a connected finite acyclic non-Dynkin quiver $Q$.
  Let $I$ denote the set of vertices of $Q$.
  Since the graded Nakayama automorphism $\mu$ of $B$ is determined only up to an inner automorphism,
  it follows from \ref{rem:graded-iso-idem} that we may assume $\mu(I)=I$.
  By \cite[Thm.\@ 4.12(2)]{minamotoStructureASGorensteinAlgebras2011}, there is an isomorphism of graded algebras
  $\phi\colon B_\mu\to\Pi(\Lambda)$.
  By \ref{rem:graded-iso-idem}, we may assume $\phi(I)=I$.
  By the definition of a graded twist,
  $\tau\defeq\phi\mu^{-1}\phi^{-1}$ defines an automorphism of $\Pi(\Lambda)$
  such that $B\cong \Pi(\Lambda)_\tau$.
  Since $\tau(I)=I$, $B$ is graded isomorphic to a twisted preprojective algebra of type $Q$.

  \textit{(2)$\Rightarrow$(1).}
  By \ref{enum:prop-preproj-CY,enum:prop-zhang-twist-CY},
  $B$ is a graded skew CY algebra of dimension 2 with Gorenstein parameter 1.

  For the remaining assertions, the definition of a graded twist gives $B_0\cong\kk Q$.
  By \ref{enum:prop-zhang-twist-noeth}, $B$ is noetherian if and only if $\Pi(\kk Q)$ is noetherian, which is equivalent to $Q$ being Euclidean by \ref{prop:preproj}.
\end{proof}

We will also use the Hilbert series of twisted preprojective algebras.

\begin{definition}[{{\cite[Sec.\@ 2.11]{beilinsonKoszulDualityPatterns1996}}}]
  Let $Q$ be a finite quiver.
  Denote its vertex set by $Q_0=\{e_0,e_1,\dotsc,e_{n-1}\}$.
  The \dfn{adjacency matrix} of $Q$
  is the matrix $C\in\uM_n(\ZZ)$ whose $(i,j)$-entry
  is the number of arrows from $e_i$ to $e_j$.

  Recall that $\kk Q_0=\bigoplus_{e\in Q_0}\kk e$.
  Let $M$ be a locally finite $\NN$-graded $\kk Q_0$-bimodule.
  The \dfn{Hilbert series} of $M$
  is the matrix-valued series $H_M(t)\in \uM_n(\ZZ[[t]])$
  with entries
  \begin{equation*}
    H_M(t)_{ij}=\sum_{m=0}^{\infty}\dim_{\kk} (e_iM_me_j)t^m.
  \end{equation*}
\end{definition}

We begin with the Hilbert series of graded preprojective algebras.
A direct computation using \ref{enum:prop-preproj-res} establishes the following formula;
see also \cite[Prop.\@ 7.1(3)]{reyesGrowthGradedTwisted2019}.

\begin{proposition}
  Let $Q$ be a connected finite acyclic non-Dynkin quiver.
  Then the Hilbert series of its graded preprojective algebra $\Pi$
  is given by
  \[H_\Pi(t)=(I_n-C-tC^{\intercal}+tI_n)^{-1},\]
  where $I_n$ is the identity matrix in $\uM_n(\ZZ)$
  and $C^{\intercal}$ is the transpose of $C$.
\end{proposition}

Let $B=\Pi_\tau$ be a $\tau$-twisted preprojective algebra of type $Q$,
and let $\sigma$ be the permutation of $\{0,\dotsc,n-1\}$
defined by $\tau(e_j)=e_{\sigma(j)}$.
By definition,
\(
  e_iB_me_j=e_i\Pi_m e_{\sigma^m(j)}.
\)
Hence, the matrix Hilbert series of $B$ is given by
\[
  H_B(t)_{ij}=\sum_{m=0}^{\infty}\dim_{\kk}(e_i\Pi_m e_{\sigma^m(j)})t^m.
\]
In particular, if $\tau(e_i)=e_i$, then $H_B(t)_{ii}=H_\Pi(t)_{ii}$,
so $e_iBe_i$ and $e_i\Pi e_i$ have the same Hilbert series.

\subsection{Noncommutative resolutions of noncommutative isolated quadric surface singularities}

In this subsection, we establish a correspondence between
noncommutative isolated quadric surface singularities
and certain twisted preprojective algebras.
These twisted preprojective algebras give noncommutative resolutions of the singularities.

We first fix orientations of the Euclidean quivers
$\widetilde{A}_1$ and $\widetilde{D}_4$ and describe their graded preprojective algebras.

\begin{example}\label{ex:A1}
  Consider the Euclidean quiver $\widetilde{A_1}$.
  \[\begin{tikzcd}
    0 && 1
    \arrow["\beta"{description}, shift left=-2,
      from=1-1, to=1-3, leftarrow]
    \arrow["\alpha"{description}, shift left=2,
      from=1-1, to=1-3, leftarrow]
  \end{tikzcd}\tag{$\widetilde{A_1}$}\]
  Its preprojective algebra is
  $\Pi=\kk\widetilde{A_1}^{\raisemath{-2pt}{\sharp}}
  \div(\alpha\alpha^*+\beta\beta^*,\alpha^*\alpha+\beta^*\beta)$.
  The Hilbert series of $\Pi$ is
  \begin{equation*}
    H_{\Pi}(t)=(1-t)^{-2}
    \begin{bmatrix}
      1+t & 2t\\
      2 & 1+t
    \end{bmatrix}.
  \end{equation*}
\end{example}

\begin{example}\label{ex:D4}
  Consider the Euclidean quiver
  $\widetilde{D_4}$.
  \[\begin{tikzcd}[row sep=small]
    1 && 2 \\
    & 0 \\
    4 && 3
    \arrow["{\alpha_1}"{description}, from=2-2, to=1-1, leftarrow]
    \arrow["{\alpha_2}"{description}, from=2-2, to=1-3, leftarrow]
    \arrow["{\alpha_3}"{description}, from=2-2, to=3-3, leftarrow]
    \arrow["{\alpha_4}"{description}, from=2-2, to=3-1, leftarrow]
  \end{tikzcd}\tag{$\widetilde{D_4}$}\]
  Its preprojective algebra is
  $\Pi=\kk\widetilde{D_4}^{\raisemath{-2pt}{\sharp}}
  \div(\alpha_1\alpha_1^*,\alpha_2\alpha_2^*,
  \alpha_3\alpha_3^*,\alpha_4\alpha_4^*,
  \sum_{i=1}^4\alpha_i^*\alpha_i
  )$.
  The Hilbert series of $\Pi$ is
  \begin{equation*}
    H_{\Pi}(t)=(1-t)^{-2}(1+t)^{-1}
    \begin{bmatrix}
      (1+t)^2 & t(1+t) & t(1+t) & t(1+t) & t(1+t)\\
      1+t & 1-t+t^2 & t & t & t\\
      1+t & t & 1-t+t^2 & t & t\\
      1+t & t & t & 1-t+t^2 & t\\
      1+t & t & t & t & 1-t+t^2
    \end{bmatrix}.
  \end{equation*}
\end{example}

For $d=2$, the following lemma identifies the algebra $\Lambda$
in \ref{thm:tilting-hypersurface} with the path algebra
of one of the quivers in \ref{ex:A1,ex:D4}.

\begin{lemma}\label{cor:tilting-surface}
  Let $A$ be a noncommutative isolated quadric surface singularity.
  Then the algebra $\Lambda$ appearing in \ref{thm:tilting-hypersurface}
  is isomorphic to either $\kk\widetilde{A_1}$ or $\kk\widetilde{D_4}$.
\end{lemma}

\begin{proof}
  Keep the notation in \ref{thm:tilting-hypersurface}.
  Recall that $\Dbqgr{A}$ has a full strong exceptional sequence
  $\bigl(\pi(A),\pi(M_1),\dotsc,\pi(M_n)\bigr)$,
  and that $\Hom_{\qgr{A}}(\pi(M_i),\pi(M_j))=0$
  whenever $i\neq j$.
  It follows that $\Lambda\cong \kk Q$ for a quiver $Q$
  with vertex set $\{0, \dotsc, n\}$
  such that every arrow goes from a vertex $i\in\{1,\dotsc,n\}$ to vertex $0$.
  For $i=1,\dotsc, n$, 
  \[
    e_i\Lambda e_0\cong\Hom_{\qgr A}(\pi(A),\pi(M_i))\cong\Hom_{\gr A}(A, M_i)
    \cong (M_i)_0.
  \]
  Thus, $Q$ has $\dim_{\kk}((M_i)_0)$ arrows from vertex $i$ to vertex $0$.
  Since $C(A)$ is semisimple and $\dim_{\kk} (C(A))=4$
  by \ref{thm:quadric-iso-sing,enum:lem-CA-dim},
  $C(A)$ is isomorphic to either $M_2(\kk)$ or $\kk^4$.
  By \ref{enum:lem-CM-module-quadric-corr},
  $Q$ is the $\widetilde{A_1}$ quiver when $C(A)\cong M_2(\kk)$
  and the $\widetilde{D_4}$ quiver when $C(A)\cong \kk^4$.
\end{proof}

Every noncommutative isolated quadric surface singularity admits a noncommutative resolution
in the form of a twisted preprojective algebra.
The following proposition summarizes the relevant results.

\begin{proposition}\label{prop:NCCR}
  Let $A$ be a noncommutative isolated quadric surface singularity,
  and let $M$ denote the direct sum of all pairwise non-isomorphic indecomposable non-projective modules in $\CM*{A}$.
  Let $B=\gEnd_A(A\oplus M)$.
  Then the following statements hold.
  \begin{enumerate}[ref={\ref{prop:NCCR}(\arabic*)},beginpenalty=10000]
    \item $A\oplus M$ is a graded 1-cluster tilting module.
      Consequently, $B$ is a noncommutative resolution of $A$.
    \item $B$ is $\NN$-graded skew CY of dimension $2$
      with Gorenstein parameter $1$.
    \item\label{enum:prop-NCCR-iso}
      $B$ is graded isomorphic to a twisted preprojective algebra
      of type $\widetilde{A_1}$ or $\widetilde{D_4}$.
  \end{enumerate}
\end{proposition}

\begin{proof}
  (1) It follows directly from the definition and \ref{thm:NCCR-ct}, since $\add(A\oplus M)=\CM{A}$.

  (2) Since $A$ has Gorenstein parameter 1, it follows from the proof of \cite[Thm.\@ 3.10(1)]{ueyamaClusterTiltingModules2017} and \ref{thm:CY-equiv}.
  
  (3) It follows from (2), \ref{cor:tilting-surface}, and \ref{cor:twisted-preproj}.
\end{proof}



Conversely, let $B$ be a twisted preprojective algebra of type
$\widetilde{A_1}$ or $\widetilde{D_4}$.
We show that $e_0Be_0$ is a noncommutative isolated quadric surface singularity.
To this end, we use the following lemma,
which follows directly from the definition of a graded twist.

\begin{lemma}\label{lem:twist-restriction}
  Let $B=\Pi_\tau$ be a $\tau$-twisted preprojective algebra.
  If $\tau(e_i)=e_i$ for some idempotent $e_i\in\Pi_0$,
  then $\tau$ restricts to an automorphism of $e_i \Pi e_i$,
  and $e_i B e_i$ is graded isomorphic
  to the twisted algebra $(e_i \Pi e_i)_\tau$.
\end{lemma}



\begin{proposition}\label{prop:eBe-quadric}
  Let $B$ be a $\tau$-twisted preprojective algebra of type
  $\widetilde{A_1}$ or $\widetilde{D_4}$.
  Then $A\defeq e_0Be_0$ is a
  noncommutative isolated quadric surface singularity.
\end{proposition}

\begin{proof}
  First, consider the type $\widetilde{A_1}$.
  Let $\Lambda=\kk\widetilde{A_1}$, let $\Pi=\Pi(\Lambda)$,
  and let $A'\defeq e_0\Pi e_0$.
  Let $S'=\kk[x,y,z]$ and let $f=xz+y^2\in S'$.
  There is a surjection $p\colon \kk[x,y,z]\div(xz+y^2)\to A'$
  given by $x\mapsto \beta^*\alpha, y\mapsto \alpha^*\alpha,
  z\mapsto \alpha^*\beta$.
  Since both algebras have the same Hilbert series
  $(1-t)^{-2}(1+t)$ by \ref{enum:hilbert-quadric,ex:A1},
  $p$ is an isomorphism.
  Therefore, $A'$ is a (commutative) isolated quadric surface singularity.
  Since $e_0$ is the unique $e_j$ such that
  $e_i\Pi_0e_j\neq 0$ for all $i$,
  it follows that $\tau(e_0)=e_0$.
  By \ref{lem:twist-restriction}, the given automorphism $\tau$ of $\Pi$ restricts to
  an automorphism $\tau_{A'}$ of $A'$.
  Since $\tau_{A'}$ restricts to a $\kk$-linear automorphism of $A_1'$,
  it lifts to a graded automorphism $\tau_{S'}$ of $S'$
  such that the following diagram commutes,
  \begin{equation*}
    \begin{tikzcd}
      S' & S'\\
      A' & A'
      \arrow[from=1-1,]{1-2}[]{\tau_{S'}}
      \arrow[from=2-1,]{2-2}[]{\tau_{A'}}
      \arrow[from=1-1,]{2-1}[]{\pi}
      \arrow[from=1-2,]{2-2}[]{\pi}
    \end{tikzcd}
  \end{equation*}
  where $\pi$ is the composite of the canonical projection $S'\to S'\div(f)$
  with $p$.
  It follows that $\tau_{S'}(f)=\lambda f$ for some $\lambda\in \kk^\times$.
  Therefore, $A\cong A_{\tau_{A'}}'\cong S'_{\tau_{S'}}\div(f_{\tau_{S'}})$
  by \ref{lem:twist-restriction,lem:twisted-alg},
  and $f_{\tau_{S'}}$ is normal regular by \ref{lem:twist-normal-regular}.
  By \ref{enum:prop-zhang-twist-noeth,enum:prop-zhang-twist-Koszul},
  $A$ is a noncommutative isolated quadric surface singularity.

  Next, consider the type $\widetilde{D_4}$.
  Let $\Lambda=\kk\widetilde{D_4}$, let $\Pi=\Pi(\Lambda)$,
  and let $A'\defeq e_0\Pi e_0$.
  Denote by $\kk_{-1}[x,y,z]$
  the $(-1)$-skew polynomial algebra
  $\kk\langle x,y,z\rangle\div(xy+yx,xz+zx,yz+zy)$.
  The element $f=x^2+y^2+z^2$ is central and regular in $\kk_{-1}[x,y,z]$.
  There is a surjection
  \begin{gather*}
    p\colon \kk_{-1}[x,y,z]\div(f)\to A',\\
    x\mapsto \alpha_1^*\alpha_1+\alpha_2^*\alpha_2,\;
    y\mapsto \alpha_1^*\alpha_1+\alpha_3^*\alpha_3,\;
    z\mapsto \alpha_2^*\alpha_2+\alpha_3^*\alpha_3.
  \end{gather*}
  By \ref{enum:hilbert-quadric,ex:D4}, both algebras have the same Hilbert series $(1-t)^{-2}(1+t)$.
  Thus, $p$ is an isomorphism
  and $A'$ is a noncommutative isolated quadric surface singularity
  by \cite[Thm.\@ 3.10]{ueyamaKnorrerPeriodicityQuadric2019}
  and \ref{thm:quadric-iso-sing}.

  We claim that $\tau$ permutes the four elements $\alpha_i^*\alpha_i$
  up to multiplication by a common scalar.
  Note that
  \begin{equation*}
    \dim_{\kk}(e_i\Pi_0e_0)=\dim_{\kk}(e_0\Pi_1e_i)=1
    \text{ and }
    \dim_{\kk}(e_i\Pi_0e_j)=0
  \end{equation*}
  for $i\neq j\in\{1,2,3,4\}$.
  Since $\tau$ restricts to $\kk$-linear bijections
  $e_i\Pi_0e_0\cong \tau(e_i)\Pi_0\tau(e_0)$,
  it follows that $\tau(e_0)=e_0$
  and that there exists a permutation $\sigma\in S_4$
  such that $\tau(e_i)=e_{\sigma(i)}$ for $i\in\{1,2,3,4\}$.
  Moreover, there exist scalars $\lambda_i, \lambda_i'\in \kk^\times$
  such that $\tau(\alpha_i)=\lambda_i\alpha_{\sigma(i)}$
  and $\tau(\alpha_i^*)=\lambda_i' \alpha_{\sigma(i)}^*$.
  Since $\tau$ lifts to an automorphism $\Tau$ of
  the path algebra $\kk\widetilde{D_4}^{\raisemath{-2pt}{\sharp}}$
  that preserves the defining ideal of $\Pi$,
  there exists a scalar $\lambda\in \kk^\times$ such that
  \begin{equation*}
    \lambda\sum_{i=1}^4 \alpha_i^* \alpha_i
    =\Tau\Biggl(\sum_{i=1}^4 \alpha_i^* \alpha_i\Biggr)
    =\sum_{i=1}^4 \lambda_i \lambda_i' \alpha_{\sigma(i)}^* \alpha_{\sigma(i)}.
  \end{equation*}
  It follows that $\lambda_i'=\lambda \lambda_i^{-1}$ for $i\in\{1,2,3,4\}$.

  Since $\sum_{i=1}^4\alpha_i^*\alpha_i=0$,
  the restriction $\tau_{A'}\defeq\tau|_{A'}$ permutes $p(x),p(y),p(z)$
  up to signs and multiplication by $\lambda$.
  Thus, $\tau_{A'}$ lifts to a graded automorphism of $\kk_{-1}[x,y,z]$ that
  sends $f$ to $\lambda^2f$.
  The argument used in the type $\widetilde{A_1}$ case also implies that
  $A$ is a noncommutative isolated quadric surface singularity.
\end{proof}

It remains to show that
each algebra $B$ in \ref{prop:eBe-quadric} is a noncommutative resolution
of $e_0Be_0$.
To this end, we use the Morita-like theory for graded quotient categories $\qgrOp$
developed in \cite{liNoncommutativeResolutionsASGorenstein2026}
and further studied in \cite{liNoncommutativeResolutionsNoncommutative2026}.

\begin{theorem}\label{thm:eBe-morita}
  Let $B$ be a noetherian $\NN$-graded twisted CY algebra
  of dimension $d\geq 2$.
  Suppose that $e\in B_0$ is an idempotent such that
  $\dim_{\kk} (B\div (e))$ is finite,
  and let $A=eBe$.
  Then $A$ is noetherian, and the following statements hold.
  \begin{enumerate}
    \item The $\Ao$-module $eB$ and the $A$-module $Be$
      are finitely generated.
    \item There are mutually quasi-inverse functors
      $\blank\otimes_B Be\colon \qgr{B}\simeq\qgr{A}
      \reflectbox{$\colon$}\blank\otimes_A eB$.
    \item The $A$-module $Be$ has depth $d$.
    \item There is an isomorphism of graded algebras
      $B\cong\gEnd_A(Be)$.
  \end{enumerate}
\end{theorem}

\begin{proof}
  (1) It follows from
  \cite[Lem.\@ 5.1]{liNoncommutativeResolutionsASGorenstein2026}.

  (2) By (1), the result follows from
  \cite[Lem.\@ 2.7]{moriAmpleGroupAction2016}.

  (3) It follows from
  \cite[Prop.\@ 5.2]{liNoncommutativeResolutionsNoncommutative2026}.

  (4) It follows from (3) and
  \cite[Thm.\@ 4.5(3)]{liNoncommutativeResolutionsASGorenstein2026}.
\end{proof}

Using \ref{thm:eBe-morita}, we obtain the following.

\begin{theorem}\label{thm:classif-corresp}
  Let $A$ be an $\NN$-graded algebra. The following conditions are equivalent.
  \begin{enumerate}
    \item $A$ is a noncommutative isolated quadric surface singularity.
    \item $A\cong e_0Be_0$ as graded algebras for a twisted preprojective algebra $B$
      of type $\widetilde{A_1}$ or $\widetilde{D_4}$, where $B$ is unique up to graded isomorphism.
  \end{enumerate}
  
  Consequently, there is a one-to-one correspondence
  between the graded isomorphism classes of
  noncommutative isolated quadric surface singularities
  and the graded isomorphism classes of
  twisted preprojective algebras of type $\widetilde{A_1}$ or $\widetilde{D_4}$.
\end{theorem}

\begin{proof}
  \textit{(1)$\Rightarrow$(2).}
  Let $M_A$ be as in \ref{prop:NCCR}, and let $B=\gEnd_A(A\oplus M)$.
  By \ref{enum:prop-NCCR-iso}, we may identify $B$ with a twisted preprojective algebra
  of type $\widetilde{A_1}$ or $\widetilde{D_4}$.
  By \ref{rem:graded-iso-idem}, we may assume that the $e_i\in B_0$
  are the projection maps from $A\oplus M$ onto indecomposable direct summands.
  Since $e_0$ is the unique $e_j$ such that
  $e_iB_0e_j\neq 0$ for all $i$,
  it corresponds to the summand $A$ of $A\oplus M$.
  Hence $A\cong e_0Be_0$ as graded algebras.

  To prove uniqueness, suppose that
  $A\cong e_0Be_0\cong e_0'B'e_0'$ as graded algebras,
  where $B$ and $B'$ are twisted preprojective algebras
  of type $\widetilde{A_1}$ or $\widetilde{D_4}$.
  Identify $A$ with $e_0Be_0$.
  Fix a graded isomorphism from $A$ to $e_0'B'e_0'$
  and use it to regard $B'e_0'$ as a graded $A$-module.

  By \ref{cor:twisted-preproj},
  $B$ is a noetherian graded skew CY algebra of dimension 2.
  Since $\dim_{\kk} (B\div (e_0))$ is finite,
  $Be_0$ is a Cohen--Macaulay $A$-module by \ref{thm:eBe-morita}(3).
  Moreover, by \ref{thm:eBe-morita}(4), $B\cong\gEnd_A(Be_0)$.
  The module $Be_0\in\Mod{A}$ is a generator of $\Mod{A}$ because it contains $A=e_0Be_0$ as a direct summand.
  Thus, $B$ is a noncommutative resolution of $A$.
  By \ref{thm:NCCR-ct}, $Be_0$ is an additive generator of $\CM{A}$.

  We next show that $Be_0$ is basic and generated in degree 0.
  Let $Q$ denote the quiver of $B_0$, and
  put $N_i=e_iBe_0$ for $i\in Q_0$.
  Since $\End_{\gr{A}}(Be_0)\cong B_0=\kk Q$ is basic,
  the modules $N_i$ are pairwise non-isomorphic indecomposable
  Cohen--Macaulay modules.
  By \ref{lem:CM-shift}, every $N_i$ is generated
  in a single degree.
  Since $B$ is $\NN$-graded, $(N_i)_m=0$ for $m<0$. Moreover, $(N_i)_0=e_iB_0e_0=e_i(\kk Q)e_0\neq 0$ for every $i$. Thus, each $N_i$ is generated in degree 0.

  The same argument applies to $B'e_0'$.
  Thus $Be_0$ and $B'e_0'$ are basic additive generators of $\CM{A}$
  whose indecomposable summands are generated in degree 0. Consequently, they are isomorphic as graded modules.
  In particular, $B$ and $B'$ have the same quiver type and are graded isomorphic.

  \textit{(2)$\Rightarrow$(1).}
  This follows from \ref{prop:eBe-quadric}.
\end{proof}

\subsection{Classification of twisted preprojective algebras of types \texorpdfstring{$\widetilde{A}_1$}{extended A1} and \texorpdfstring{$\widetilde{D}_4$}{extended D4}}

By \ref{thm:classif-corresp}, noncommutative isolated quadric surface singularities
correspond to twisted preprojective algebras of type
$\widetilde{A}_1$ or $\widetilde{D}_4$.
We now classify these algebras up to graded isomorphism.

For a nonsingular matrix $M$, we define
$M^{-\intercal}\defeq (M^{-1})^\intercal=(M^\intercal)^{-1}$.

\begin{proposition}\label{prop:B-classif-A1}
  Let $B$ be a $\tau$-twisted preprojective algebra of type $\widetilde{A_1}$.
  \begin{enumerate}[beginpenalty=10000,ref={\ref{prop:B-classif-A1}(\arabic*)}]
    \item
      $\tau$ is uniquely determined by a matrix $M\in \GL_2(\kk)$
      and a scalar $\lambda\in \kk^\times$ such that
      \begin{equation}\label{eq:tau-A1}
        \begin{bmatrix}\tau(e_0)\\\tau(e_1)\end{bmatrix}
        =\begin{bmatrix}e_0\\e_1\end{bmatrix},\quad
        \begin{bmatrix}\tau(\alpha)\\\tau(\beta)\end{bmatrix}
        =M\begin{bmatrix}\alpha\\\beta\end{bmatrix},\quad
        \begin{bmatrix}\tau(\alpha^*)\\\tau(\beta^*)\end{bmatrix}
        =N\begin{bmatrix}\alpha^*\\\beta^*\end{bmatrix},
      \end{equation}
      where $N=\lambda M^{-\intercal}$.

    \item\label{enum:prop-B-classif-A1-2}
      $B$ is isomorphic, as a graded algebra, to $\symsf{A}(M^{-1})$, which is
      defined by the weighted quiver
      $\widetilde{A_1}^{\raisemath{-2pt}{\sharp}}$
      with the two relations
      \begin{equation*}
        \begin{bmatrix}\alpha&\beta\end{bmatrix}\begin{bmatrix}\alpha^*\\\beta^*\end{bmatrix},
        \quad
        \begin{bmatrix}\alpha^*&\beta^*\end{bmatrix}M^{-1}\begin{bmatrix}\alpha\\\beta\end{bmatrix}.
      \end{equation*}
      In particular, $B$ is independent of the choice of $\lambda$.
    \item
      Up to graded isomorphism,
      the $\tau$-twisted preprojective algebras of
      type $\widetilde{A_1}$ are represented by the graded algebras
      $\symsf{A}_1(p)\defeq \symsf{A}(M_p^{-1})$
      for $p\in \kk^\times$ and $\symsf{A}_2\defeq \symsf{A}(M^{-1})$,
      where
      \begin{equation*}
        M_p^{-1}=\begin{bmatrix}1\\&p\end{bmatrix}\text{ and }
        M^{-1}=\begin{bmatrix}1&1\\&1\end{bmatrix}.
      \end{equation*}
      These algebras are pairwise non-isomorphic as graded algebras, except that
      $\symsf{A}_1(p)\cong\symsf{A}_1(p^{-1})$ for $p\in \kk^\times$.
      %
  \end{enumerate}
\end{proposition}

\begin{proof}
  (1) Let $\Pi$ denote the preprojective algebra of $\kk\widetilde{A_1}$.
  Note that
  \begin{equation*}
    \dim_{\kk}(e_1\Pi_0e_0)=\dim_{\kk}(e_0\Pi_1e_1)=2
    \text{ and }
    \dim_{\kk}(e_0\Pi_0e_1)=0.
  \end{equation*}
  Since $\tau$ restricts to $\kk$-linear bijections
  \begin{equation*}
    \kk\alpha\oplus \kk\beta=e_1\Pi_0e_0\cong \tau(e_1)\Pi_0\tau(e_0)
    \text{ and }
    \kk\alpha^*\oplus \kk\beta^*=e_0\Pi_1e_1\cong\tau(e_0)\Pi_1\tau(e_1),
  \end{equation*}
  it follows that $\tau$ has the form given in \ref{eq:tau-A1},
  where $M,N\in \GL_2(\kk)$.
  It remains to show that $N=\lambda M^{-\intercal}$ for some $\lambda\in \kk^\times$.
  Note that $\tau$ lifts to an automorphism $\Tau$ of
  $\kk\widetilde{A_1}^{\raisemath{-2pt}{\sharp}}$,
  which preserves the defining ideal.
  Thus, there exists a scalar $\lambda\in \kk^\times$ such that
  \begin{equation*}
    \lambda(\alpha\alpha^*+\beta\beta^*)=\Tau(\alpha\alpha^*+\beta\beta^*)
    =\begin{bmatrix}\alpha&\beta\end{bmatrix}M^\intercal N
    \begin{bmatrix}\alpha^*\\\beta^*\end{bmatrix},
  \end{equation*}
  that is, $M^\intercal N=\lambda I_2$.
  Similarly, there exists $\lambda'\in \kk^\times$
  such that $N^\intercal M=\lambda' I_2$.
  Hence $\lambda=\lambda'$ and $N=\lambda M^{-\intercal}$.

  (2) It follows from a computation using
  \ref{ex:A1,lem:twisted-alg}.

  (3) Let $B'$ be a $\tau'$-twisted preprojective algebra of type $\widetilde{A_1}$,
  and suppose that $B$ and $B'$ are graded isomorphic.
  By (2), we may assume that $B=\symsf{A}(M^{-1})$ and $B'=\symsf{A}((M')^{-1})$.
  By \ref{rem:graded-iso-idem}, there is a graded isomorphism
  $\phi\colon B\to B'$ that permutes the vertices of $\widetilde{A_1}$.
  Then by the same argument as in (1),
  $\phi$ is determined by matrices $U,V\in \GL_2(\kk)$ such that
  \begin{equation}\label{eq:phi-A1}
    \begin{bmatrix}\phi(e_0)\\\phi(e_1)\end{bmatrix}
    =\begin{bmatrix}e_0'\\e_1'\end{bmatrix},\quad
    \begin{bmatrix}\phi(\alpha)\\\phi(\beta)\end{bmatrix}
    =U\begin{bmatrix}\alpha'\\\beta'\end{bmatrix},\quad
    \begin{bmatrix}\phi(\alpha^*)\\\phi(\beta^*)\end{bmatrix}
    =V\begin{bmatrix}\alpha^{*\prime}\\\beta^{*\prime}\end{bmatrix}.
  \end{equation}
  Note that $\phi$ lifts to an automorphism $\Phi$
  of $\kk\widetilde{A_1}^{\raisemath{-2pt}{\sharp}}$,
  which maps the defining ideal of $B$ onto that of $B'$.
  Thus, there exist scalars $\gamma,\delta\in \kk^\times$ such that
  \begingroup\allowdisplaybreaks
  \begin{align*}
    \gamma(\alpha'\alpha^{*\prime}+\beta'\beta^{*\prime})
    =\Phi(\alpha\alpha^*+\beta\beta^*)
    &=\begin{bmatrix}\alpha'&\beta'\end{bmatrix}U^\intercal V
    \begin{bmatrix}\alpha^{*\prime}\\\beta^{*\prime}\end{bmatrix},\\
    \delta\begin{bmatrix}\alpha^{*\prime}&\beta^{*\prime}\end{bmatrix}(M')^{-1}
    \begin{bmatrix}\alpha'\\\beta'\end{bmatrix}
    =\Phi\Biggl(\begin{bmatrix}\alpha^*&\beta^*\end{bmatrix}M^{-1}
    \begin{bmatrix}\alpha\\\beta\end{bmatrix}\Biggr)
    &=\begin{bmatrix}\alpha^{*\prime}&\beta^{*\prime}\end{bmatrix}V^\intercal M^{-1}U
    \begin{bmatrix}\alpha^{\prime}\\\beta^{\prime}\end{bmatrix}.
  \end{align*}
  \endgroup
  It follows that $U^\intercal V=\gamma I_2$
  and $V^\intercal M^{-1}U=\delta(M')^{-1}$.
  Therefore,
  \begin{equation*}
    (M')^{-1}=\gamma\delta^{-1} U^{-1} M^{-1} U.
  \end{equation*}

  Conversely, suppose that $(M')^{-1}=\lambda U^{-1}M^{-1}U$
  for some $\lambda\in\kk^\times$ and $U\in\GL_2(\kk)$.
  Setting $V=U^{-\intercal}$, the assignments in \ref{eq:phi-A1}
  define a graded isomorphism from $B$ to $B'$.

  Thus, $B$ and $B'$ are graded isomorphic
  if and only if $M^{-1}$ and $(M')^{-1}$ are similar up to a scalar.
  Equivalently, $M^{-1}$ and $(M')^{-1}$ are conjugate
  in $\PGL_2(\kk)$.
  The Jordan normal form then yields the desired result.
\end{proof}

\begin{proposition}\label{prop:B-classif-D4}
  Let $B$ be a $\tau$-twisted preprojective algebra of type $\widetilde{D_4}$.
  \begin{enumerate}[beginpenalty=10000,ref={\ref{prop:B-classif-D4}(\arabic*)}]
    \item\label{enum:prop-B-classif-D4-tau}
      $\tau$ is uniquely determined by a permutation
      $\sigma\in S_4=\Aut{\{1,2,3,4\}}$
      and a 5-tuple of scalars
      $(\lambda_1,\lambda_2,\lambda_3,\lambda_4,\lambda)
      \in(\kk^\times)^{\oplus 5}$
      such that
      \begin{equation*}
        \tau(e_0)=e_0,\quad \tau(e_i)=e_{\sigma(i)},
        \quad
        \tau(\alpha_i)=\lambda_i \alpha_{\sigma(i)},\quad
        \tau(\alpha_i^*)=\lambda \lambda_i^{-1} \alpha_{\sigma(i)}^*
      \end{equation*}
      for $i\in \{1,2,3,4\}$.
    \item
      $B$ is isomorphic, as a graded algebra, to
      $\symsf{D}(\sigma^{-1})$, defined by the weighted quiver
      $\widetilde{D_4}^{\raisemath{-2pt}{\sharp}}$
      and the five relations
      \begin{equation*}
        \sum_{i=1}^4 \alpha_i^* \alpha_i,
        \quad \alpha_i \alpha_{\sigma^{-1}(i)}^*
        \text{ for }i\in \{1,2,3,4\}.
      \end{equation*}
      In particular, $B$ is independent of the choice of $(\lambda_1,\lambda_2,\lambda_3,\lambda_4,\lambda)$.
    \item
      Up to graded isomorphism,
      the $\tau$-twisted preprojective algebras of
      type $\widetilde{D_4}$ are represented by the five graded algebras
      \begin{equation*}
        \symsf{D}_1\defeq \symsf{D}(\id),
        \symsf{D}_2\defeq \symsf{D}\bigl((12)\bigr),
        \symsf{D}_3\defeq \symsf{D}\bigl((12)(34)\bigr),
        \symsf{D}_4\defeq \symsf{D}\bigl((123)\bigr),
        \symsf{D}_5\defeq \symsf{D}\bigl((1234)\bigr).
      \end{equation*}
  \end{enumerate}
\end{proposition}

\begin{proof}
  (1) It follows from the computation in the type $\widetilde{D_4}$ case
  of the proof of \ref{prop:eBe-quadric}.

  (2) We sketch the proof.
  By \ref{ex:D4,lem:twisted-alg},
  $B$ is graded isomorphic to the algebra defined by
  the following weighted quiver
  \[\begin{tikzcd}
    1 && 2 \\
    & 0 \\
    4 && 3
    \arrow["{\alpha_{\sigma(1)}^*}"{description},
      curve={height=3pt}, shift right=1, from=1-1, to=2-2, leftarrow]
    \arrow["{\alpha_{\sigma(2)}^*}"{description},
      curve={height=3pt}, shift right=1, from=1-3, to=2-2, leftarrow]
    \arrow["{\alpha_{\sigma(3)}^*}"{description},
      curve={height=3pt}, shift right=1, from=3-3, to=2-2, leftarrow]
    \arrow["{\alpha_{\sigma(4)}^*}"{description},
      curve={height=3pt}, shift right=1, from=3-1, to=2-2, leftarrow]
    \arrow["{\alpha_1}"{description}, curve={height=3pt}, shift right=1,
      from=2-2, to=1-1, leftarrow]
    \arrow["{\alpha_2}"{description}, curve={height=3pt}, shift right=1,
      from=2-2, to=1-3, leftarrow]
    \arrow["{\alpha_3}"{description}, curve={height=3pt}, shift right=1,
      from=2-2, to=3-3, leftarrow]
    \arrow["{\alpha_4}"{description}, curve={height=3pt}, shift right=1,
      from=2-2, to=3-1, leftarrow]
  \end{tikzcd}\qquad
  \begin{aligned}
    w(\alpha_i)&=0,\; i\in\{1,2,3,4\};\\
    w(\alpha_i^*)&=1,\; i\in\{1,2,3,4\},
  \end{aligned}
  \]
  with the five relations
  \begin{equation*}
    \sum_{i=1}^4 \lambda_{\sigma^{-1}(i)}^{-1}\alpha_i^* \alpha_{\sigma^{-1}(i)},
    \quad \alpha_i \alpha_i^*
    \text{ for }i\in \{1,2,3,4\}.
  \end{equation*}
  Replacing $\alpha_i^*$
  by $\lambda_{\sigma^{-1}(i)}\alpha_{\sigma^{-1}(i)}^*$
  for $i\in\{1,2,3,4\}$,
  we obtain $B\cong\symsf{D}(\sigma^{-1})$.

  (3) Let $B'$ be a $\tau'$-twisted preprojective algebra of type $\widetilde{D_4}$,
  and suppose that $B$ and $B'$ are graded isomorphic.
  By (2), we may assume that $B=\symsf{D}(\sigma^{-1})$ and $B'=\symsf{D}(\sigma^{\prime-1})$.
  By \ref{rem:graded-iso-idem}, there is a graded isomorphism
  $\phi\colon B\to B'$ that permutes the vertices of $\widetilde{D_4}$.
  Similar to (1),
  $\phi$ is determined by a permutation $\rho\in S_4$
  and two 4-tuples of scalars $\theta,\theta'\in(\kk^\times)^{\oplus 4}$
  such that
  \begin{equation*}
    \phi(e_0)=e_0',\; \phi(e_i)=e_{\rho(i)}',\;
    \phi(\alpha_i)=\theta_i \alpha_{\rho(i)}',\;
    \phi(\alpha_i^*)=\theta_i' {\alpha_{\rho(i)}^{*\prime}}.
  \end{equation*}
  Comparing the relations for paths of length 2 starting at vertices $1,\dotsc,4$,
  we obtain
  \[
  \{\alpha_{\rho(i)}'\alpha_{\rho\sigma^{-1}(i)}^{*\prime}
    \mid i=1,\dotsc,4\}
  =
  \{\alpha_i'\alpha_{\sigma^{\prime-1}(i)}^{*\prime}
    \mid i=1,\dotsc,4\}.
  \]
  It follows that $\sigma^{\prime-1}=\rho\sigma^{-1}\rho^{-1}$.
  Conversely, such a $\rho$ defines a graded isomorphism from $B$ to $B'$,
  as can be verified by taking $\theta_i=\theta_i'=1$ for every $i\in\{1,2,3,4\}$.
  The result then follows from the classification of
  conjugacy classes of $S_4$.
\end{proof}

The preceding two propositions lead to our first classification theorem.

\begin{theorem}\label{thm:B-classif}
  Let $B$ be a $\tau$-twisted preprojective algebra of type
  $\widetilde{A_1}$ or $\widetilde{D_4}$.
  Then $B$ is graded isomorphic to one of the algebras listed in
  \ref{tab:B-classif}.
  The listed algebras are pairwise non-isomorphic as graded algebras,
  except that $\symsf{A}_1(p)\cong\symsf{A}_1(p^{-1})$
  for $p\in\kk^\times$.
\end{theorem}

\begin{table}[htb]
  \caption{The classification (\ref{thm:B-classif})
    of twisted preprojective algebras
    $B=\kk Q^\sharp\div I$
    of type $\widetilde{A_1}$ or $\widetilde{D_4}$ up to graded isomorphism
    consists of one family and six isolated cases.
    In the family $\symsf{A}_1(p)$, the parameter $p$
    is considered up to inversion.}
  \label{tab:B-classif}
  \centering
  \begin{tabular}{c@{\hskip 2em}ccc}
    \toprule
    $Q^\sharp$ & $B$ & $I$\\
    \midrule
    \multirow{2}{*}{$\widetilde{A_1}^{\raisemath{-2pt}{\sharp}}$}
    & $\symsf{A}_1(p)$ & $(\alpha\alpha^*+\beta\beta^*,\; \alpha^*\alpha+p\beta^*\beta)$
    & $p\in \kk^\times$\\
    & $\symsf{A}_2$ & $(\alpha\alpha^*+\beta\beta^*,\; \alpha^*\alpha+\alpha^*\beta+\beta^*\beta)$\\
    \cmidrule(r{1.5em}){1-1}
    \multirow{5}{*}{$\widetilde{D_4}^{\raisemath{-2pt}{\sharp}}$}
    & $\symsf{D}_1$
    & \multicolumn{2}{@{}c}{$(\alpha_1\alpha_1^*,\; \alpha_2\alpha_2^*,\; \alpha_3\alpha_3^*,\; \alpha_4\alpha_4^*,\;
    \sum_{i=1}^4 \alpha_i^*\alpha_i)$}\\
    & $\symsf{D}_2$
    & \multicolumn{2}{@{}c}{$(\alpha_1\alpha_2^*,\; \alpha_2\alpha_1^*,\; \alpha_3\alpha_3^*,\; \alpha_4\alpha_4^*,\;
    \sum_{i=1}^4 \alpha_i^*\alpha_i)$}\\
    & $\symsf{D}_3$
    & \multicolumn{2}{@{}c}{$(\alpha_1\alpha_2^*,\; \alpha_2\alpha_1^*,\; \alpha_3\alpha_4^*,\; \alpha_4\alpha_3^*,\;
    \sum_{i=1}^4 \alpha_i^*\alpha_i)$}\\
    & $\symsf{D}_4$
    & \multicolumn{2}{@{}c}{$(\alpha_1\alpha_2^*,\; \alpha_2\alpha_3^*,\; \alpha_3\alpha_1^*,\; \alpha_4\alpha_4^*,\;
    \sum_{i=1}^4 \alpha_i^*\alpha_i)$}\\
    & $\symsf{D}_5$
    & \multicolumn{2}{@{}c}{$(\alpha_1\alpha_2^*,\; \alpha_2\alpha_3^*,\; \alpha_3\alpha_4^*,\; \alpha_4\alpha_1^*,\;
    \sum_{i=1}^4 \alpha_i^*\alpha_i)$}\\
    \bottomrule
  \end{tabular}
\end{table}

\subsection{Classification of noncommutative isolated quadric surface singularities}
\label{subsec:classif-A}

Using \ref{thm:classif-corresp}, we obtain the classification of
noncommutative isolated quadric surface singularities
by computing the corner algebras $e_0Be_0$
for the algebras $B$ in \ref{thm:B-classif}.

\begin{theorem}\label{thm:A-classif}
  There is a one-to-one correspondence between the graded isomorphism classes of noncommutative isolated quadric surface singularities
  and the union of the sets of conjugacy classes in $\PGL_2(\kk)$ and $S_4$.
  Each such singularity is graded isomorphic to one of the algebras listed in
  \ref{tab:A-classif}.
  The listed algebras are pairwise non-isomorphic as graded algebras,
  except that $\symbb{A}_1(p)\cong\symbb{A}_1(p^{-1})$
  for $p\in\kk^\times$.
\end{theorem}

\begin{table}[htb]
  \caption{The classification (\ref{thm:A-classif})
    describes noncommutative isolated quadric surface
    singularities $A=S/(f)$ up to graded isomorphism.
    The quantum polynomial algebra $S$ has dimension 3 and is of the form
    $\kk\langle x,y,z\rangle\div (f_1,f_2,f_3)$.
    The type of $A$ is the quiver $Q$ of the path algebra $\Lambda$
    appearing in \ref{cor:tilting-surface}.
    The classification consists of one family and six isolated cases.
    In the family $\symbb{A}_1(p)$, the parameter $p$
    is considered up to inversion.}
  \label{tab:A-classif}
  \centering
  \begin{tabular}{ccccc}
    \toprule
    Type & $A$ & $f_1,f_2,f_3$ & $f$\\
    \midrule
    \multirow{2}{*}{$\widetilde{A_1}$}
    & $\symbb{A}_1(p)$
    & \makecell{$pxy-yx,\;p^2xz-zx,\;pyz-zy$}
    & $pxz+y^2$ & $p\in \kk^\times$\\
    \cmidrule(r{.25em}){3-3}
    & $\symbb{A}_2$
    & \makecell{$xy+xz-yx+2y^2+yz,$\\$xz-zx+2zy+z^2,$\\$yz-zy-z^2$}
    & $xz+y^2+yz$\\
    \cmidrule{1-3}
    \multirow{5}{*}{$\widetilde{D_4}$}
    & $\symbb{D}_1$
    & \makecell{$xy+yx,\;xz+zx,\;yz+zy$}
    & $x^2+y^2+z^2$\\
    \cmidrule(r{.25em}){3-3}
    & $\symbb{D}_2$
    & \makecell{$xz+yx,\;xy+zx,\;y^2+z^2$}
    & $x^2+yz+zy$\\
    \cmidrule(r{.25em}){3-3}
    & $\symbb{D}_3$
    & \makecell{$-xy+yx,\;-xz+zx,\;yz+zy$}
    & $x^2-y^2-z^2$\\
    \cmidrule(r{.25em}){3-3}
    & $\symbb{D}_4$
    & \makecell{$x^2+yz,\;xy+z^2,\;y^2+zx$}
    & $xz+yx+zy$\\
    \cmidrule(r{.25em}){3-3}
    & $\symbb{D}_5$
    & \makecell{$-xy+yz,\;-x^2+z^2,\;yx+zy$}
    & $xz-y^2-zx$\\
    \bottomrule
  \end{tabular}
\end{table}

\begin{proof}
  By \ref{thm:classif-corresp,thm:B-classif}, it suffices to compute $e_0Be_0$
  for the algebras $B$ in \ref{tab:B-classif}.
  For $B=\Pi_\tau$, let $A'\defeq e_0\Pi e_0$
  and let $\tau_{A'}\defeq\tau|_{A'}$.
  By the proof of \ref{prop:eBe-quadric},
  $e_0Be_0\cong A'_{\tau_{A'}}$.
  The presentations in \ref{tab:A-classif} follow from \ref{lem:twisted-alg}.
\end{proof}

Let $A=S\div(f)$ and $A'=S'\div(f')$ be noncommutative quadric hypersurface singularities that are isomorphic as graded algebras. However, $f\in S$ may be noncentral while $f'\in S'$ is central. The following proposition identifies the algebras in \ref{tab:A-classif} that are isomorphic, as graded algebras, to noncommutative central quadric surface singularities.

\begin{proposition}\label{prop:A-classif-central}
  Keep the notation in \ref{tab:A-classif}.
  \begin{enumerate}
    \item
      Let $A=S\div (f)\in
      \{\mathbb{A}_1(\pm 1),\mathbb{D}_1, \mathbb{D}_2, \mathbb{D}_3\}$.
      Then $f\in S$ is central,
      and hence $A$ is a noncommutative central quadric surface singularity.
    \item Let $A=S\div (f)\in\{\mathbb{A}_1(p),\mathbb{A}_2\}$
      for $p\in \kk^\times{\setminus}\{\pm 1\}$.
      Although $f\in S$ is not central, $A$ is graded isomorphic to
      a noncommutative central quadric surface singularity.
    \item $\mathbb{D}_4$ and $\mathbb{D}_5$ are not graded isomorphic to any noncommutative central quadric surface singularity.
  \end{enumerate}
\end{proposition}

\begin{proof}
  The key to computing all degree $2$ central elements of a Koszul algebra $S$ is to apply Bergman's diamond lemma to obtain a $\kk$-basis of $S_3$. For explanations and examples, see \cite[p.\@ 29]{huClassificationNoncommutativeCentral2026}.

  (1) It follows by direct computation.

  (2) Let $A=S\div(f)=\mathbb{A}_1(p)$.
  Then $fx=p^2xf$, $fy=yf$, and $fz=p^{-2}zf$.
  Thus, $f\in S$ is not central when $p\neq \pm 1$.
  Consider $A'= S'\div(f')$, where
  \[S'\defeq \kk\langle X,Y,Z\rangle\div(pXY-YX,pYZ-ZY,pXZ+(-p^2+1)Y^2-pZX),\]
  and $f'=pXZ+Y^2$.
  The algebra $S'$ is classified as type $\symup{S_1'}$ in
  \cite[Thm.\@ 3.1]{itabaDefiningRelations3dimensional2021}. 
  The argument in \cite[Thm.\@ 3.3]{itabaASREGULARITYGEOMETRICALGEBRAS2022} also shows that
  $S'$ is a quantum polynomial algebra of dimension 3 over an arbitrary field.
  The element $f'\in S'$ is central and regular.
  The map $(x,y,z)\mapsto (X,Y,Z)$ induces an isomorphism
  of graded algebras from $A$ to $A'$.

  Let $A=S\div(f)=\mathbb{A}_2$.
  A direct computation gives
  \[
    \begin{bmatrix}fx\\fy\\fz\end{bmatrix}=
    \begin{bmatrix}
      1 & 4 & -4\\ & 1 & -2\\ & & 1
    \end{bmatrix}
    \begin{bmatrix}xf\\yf\\zf\end{bmatrix},
  \]
  so $f\in S$ is not central.
  Consider $A'=S'\div(f')$, where
  \[S'\defeq \kk\langle X,Y,Z\rangle\div(XY-YX-XZ-YZ,XZ-ZX+2YZ-Z^2,YZ-ZY-Z^2),\]
  and $f'=XY-YX+Y^2$.
  The algebra $S'$ is classified as type $\symup{T'}$ in
  \cite[Thm.\@ 3.1]{itabaDefiningRelations3dimensional2021}.
  The argument in \cite[Thm.\@ 3.3]{itabaASREGULARITYGEOMETRICALGEBRAS2022} also shows that
  $S'$ is a quantum polynomial algebra of dimension 3 over an arbitrary field.
  The element $f'\in S'$ is central and regular.
  Again, $(x,y,z)\mapsto (X,Y,Z)$ induces an isomorphism
  of graded algebras between $A$ and $A'$.

  (3) Let $A=\mathbb{D}_4$.
  Its Koszul dual is
  \[A^!=\kk\langle X,Y,Z\rangle\div(XY-ZX, -X^2+ZY, XZ-Y^2, YZ-ZX, -YX+Z^2).\]
  The degree-two part of the centre $Z(A^!)$ is spanned by
  $w\defeq X^2+XY+XZ+YX$. 
  Since $wX=wY=wZ$, the element $w$ is a zero-divisor in $A^!$.
  It follows from \ref{enum:w-central} that
  $A$ is not graded isomorphic to
  any noncommutative central quadric surface singularity.

  Let $A=\mathbb{D}_5$.
  Its Koszul dual is
  \[A^!=\kk\langle X,Y,Z\rangle\div(Y^2+ZX, YX+ZY, XZ+ZX, -XY+YZ, X^2+Z^2).\]
  The centre $Z(A^!)$ contains no nonzero element of degree 2.
  Hence, by \ref{enum:w-central},
  $A$ is not graded isomorphic to
  any noncommutative central quadric surface singularity.
\end{proof}



\subsection*{Acknowledgement}

The authors would like to thank Haigang Hu, Tiancheng Qi, and Ruipeng Zhu
for helpful discussions and comments.

\subsection*{AI statement}

Before finalizing the manuscript, the authors conducted an additional literature search using ChatGPT 6 Astra and found that, in the commutative local setting, a result similar to that in \ref{enum:high-dim-bad-depth} appears in a 2026 note by Graham J.\@ Leuschke on his website.
The authors used AI tools to refine the language, check grammar, and improve the manuscript's readability.
The mathematical results and proofs were developed by the authors, who take full responsibility for the manuscript.

  \appendix

  \myprintbibliography

@article{artinModulesRegularAlgebras1991,
  title = {Modules over Regular Algebras of Dimension 3},
  author = {Artin, Michael and Tate, John Torrence and family=Bergh, given=Michel, prefix=van den, useprefix=true},
  date = {1991-12-01},
  journaltitle = {Inventiones mathematicae},
  shortjournal = {Invent. Math.},
  volume = {106},
  number = {1},
  pages = {335--388},
  doi = {10.1007/BF01243916},
  url = {https://doi.org/10.1007/BF01243916},
  urldate = {2023-03-16},
  langid = {english},
  mrnumber = {1128218},
  shortauthor = {ATV}
}

@article{artinNoncommutativeProjectiveSchemes1994,
  title = {Noncommutative {{Projective Schemes}}},
  author = {Artin, Michael and Zhang, James Jian},
  date = {1994-12-01},
  journaltitle = {Advances in Mathematics},
  shortjournal = {Adv. Math.},
  volume = {109},
  number = {2},
  pages = {228--287},
  doi = {10.1006/aima.1994.1087},
  url = {https://www.sciencedirect.com/science/article/pii/S0001870884710875},
  urldate = {2023-03-17},
  langid = {english},
  mrnumber = {1304753}
}

@mvbook{assemElementsRepresentationTheory2006,
  title = {Elements of the {{Representation Theory}} of {{Associative Algebras}}: {{Techniques}} of {{Representation Theory}}},
  shorttitle = {Elements of the {{Representation Theory}} of {{Associative Algebras}}},
  author = {Assem, Ibrahim and Simson, Daniel and Skowronski, Andrzej},
  date = {2006},
  series = {London {{Mathematical Society Student Texts}}},
  volume = {1},
  number = {65},
  publisher = {Cambridge University Press},
  location = {Cambridge, England},
  doi = {10.1017/CBO9780511614309},
  langid = {english},
  mrnumber = {2197389},
  pagetotal = {x+458 pp.},
  shortseries = {London Math. Soc. Stud. Texts},
  volumes = {3}
}

@inproceedings{baerHomologicalPropertiesWild1986,
  title = {Homological Properties of Wild Hereditary Artin Algebras},
  booktitle = {Representation {{Theory}}},
  author = {Baer, Dagmar},
  editor = {Dlab, Vlastimil and Gabriel, Peter and Michler, Gerhard},
  date = {1986},
  series = {Lecture {{Notes}} in {{Mathematics}}},
  volume = {1},
  number = {1177},
  pages = {1--12},
  publisher = {Springer},
  location = {Berlin, Germany},
  doi = {10.1007/BFb0075255},
  eventdate = {1984-08},
  eventtitle = {The Fourth International Meeting on Representations of Algebras},
  langid = {english},
  mrnumber = {0842456},
  shortseries = {Lecture Notes in Math.},
  venue = {Ottawa}
}

@article{baerPreprojectiveAlgebraTame1987,
  title = {The Preprojective Algebra of a Tame Hereditary Artin Algebra},
  author = {Baer, Dagmar and Geigle, Werner and Lenzing, Helmut},
  date = {1987-01-01},
  journaltitle = {Communications in Algebra},
  shortjournal = {Commun. Algebra},
  volume = {15},
  number = {1--2},
  pages = {425--457},
  publisher = {Taylor \& Francis},
  doi = {10.1080/00927878708823425},
  url = {https://doi.org/10.1080/00927878708823425},
  urldate = {2025-10-14},
  mrnumber = {0876985}
}

@article{baoNoncommutativeAuslanderTheorem2018,
  title = {Noncommutative {{Auslander}} Theorem},
  author = {Bao, Yan-Hong and He, Ji-Wei and Zhang, James Jian},
  date = {2018-12},
  journaltitle = {Transactions of the American Mathematical Society},
  shortjournal = {Trans. Am. Math. Soc.},
  volume = {370},
  number = {12},
  pages = {8613--8638},
  doi = {10.1090/tran/7332},
  url = {https://www.ams.org/tran/2018-370-12/S0002-9947-2018-07332-9/},
  urldate = {2023-04-04},
  langid = {english},
  mrnumber = {3864389}
}

@article{baoPertinencyHopfActions2019,
  title = {Pertinency of {{Hopf}} Actions and Quotient Categories of {{Cohen-Macaulay}} Algebras},
  author = {Bao, Yan-Hong and He, Ji-Wei and Zhang, James Jian},
  date = {2019},
  journaltitle = {Journal of Noncommutative Geometry},
  shortjournal = {J. Noncommutative Geom.},
  volume = {13},
  number = {2},
  pages = {667--710},
  doi = {10.4171/JNCG/336},
  url = {https://mathscinet.ams.org/mathscinet-getitem?mr=3988759},
  urldate = {2023-04-28},
  mrnumber = {3988759}
}

@article{baumannCentralCurvesNoncommutative2026,
  title = {Central Curves on Noncommutative Surfaces},
  author = {Baumann, Thilo and Belmans, Pieter and Van Garderen, Okke},
  date = {2026},
  journaltitle = {Transactions of the American Mathematical Society},
  shortjournal = {Trans. Am. Math. Soc.},
  doi = {10.1090/tran/9571},
  url = {https://www.ams.org/tran/2026-379-05/S0002-9947-2025-09571-0/},
  urldate = {2026-10-30},
  langid = {english},
  mrnumber = {5060440}
}

@article{beilinsonKoszulDualityPatterns1996,
  title = {Koszul {{Duality Patterns}} in {{Representation Theory}}},
  author = {Beilinson, Alexander A. and Ginzburg, Victor A. and Soergel, Wolfgang},
  date = {1996},
  journaltitle = {Journal of the American Mathematical Society},
  shortjournal = {J. Am. Math. Soc.},
  volume = {9},
  number = {2},
  pages = {473--527},
  doi = {10.1090/S0894-0347-96-00192-0},
  url = {https://www.ams.org/jams/1996-9-02/S0894-0347-96-00192-0/},
  urldate = {2023-11-29},
  langid = {english},
  mrnumber = {1322847}
}

@article{bondalGeneratorsRepresentabilityFunctors2003,
  title = {Generators and {{Representability}} of {{Functors}} in {{Commutative}} and {{Noncommutative Geometry}}},
  author = {Bondal, Alexey I. and family=Bergh, given=Michel, prefix=van den, useprefix=true},
  date = {2003},
  journaltitle = {Moscow Mathematical Journal},
  shortjournal = {Mosc. Math. J.},
  volume = {3},
  number = {1},
  pages = {1--36},
  doi = {10.17323/1609-4514-2003-3-1-1-36},
  url = {http://www.ams.org/distribution/mmj/vol3-1-2003/bondal-vandenbergh.pdf},
  urldate = {2024-10-25},
  mrnumber = {1996800},
  shortauthor = {BV}
}

@article{bondalRepresentationAssociativeAlgebras1990,
  title = {Representation of {{Associative Algebras}} and {{Coherent Sheaves}}},
  author = {Bondal, Alexey I.},
  date = {1990-02-28},
  journaltitle = {Mathematics of the USSR-Izvestiya},
  shortjournal = {Math. USSR-Izv.},
  volume = {34},
  number = {1},
  pages = {23--42},
  publisher = {IOP Publishing},
  doi = {10.1070/IM1990v034n01ABEH000583},
  url = {https://iopscience.iop.org/article/10.1070/IM1990v034n01ABEH000583/meta},
  urldate = {2024-03-05},
  langid = {english},
  mrnumber = {0992977}
}

@article{braunFinitelyGeneratedPI1994,
  title = {Finitely {{Generated P}}.{{I}}. {{Rings}} of {{Global Dimension Two}}},
  author = {Braun, A. and Hajarnavis, C. R.},
  date = {1994-10-15},
  journaltitle = {Journal of Algebra},
  shortjournal = {J. Algebra},
  volume = {169},
  number = {2},
  pages = {587--604},
  doi = {10.1006/jabr.1994.1298},
  url = {https://www.sciencedirect.com/science/article/pii/S0021869384712981},
  urldate = {2026-10-19},
  mrnumber = {1297163}
}

@article{brownDualisingComplexesTwisted2008,
  title = {Dualising Complexes and Twisted {{Hochschild}} (Co)Homology for Noetherian {{Hopf}} Algebras},
  author = {Brown, K. A. and Zhang, J. J.},
  date = {2008-09-01},
  journaltitle = {Journal of Algebra},
  shortjournal = {J. Algebra},
  volume = {320},
  number = {5},
  pages = {1814--1850},
  doi = {10.1016/j.jalgebra.2007.03.050},
  url = {https://www.sciencedirect.com/science/article/pii/S0021869308002780},
  urldate = {2025-06-03},
  mrnumber = {2437632}
}

@article{chanMcKayCorrespondenceSemisimple2018,
  title = {{{McKay}} Correspondence for Semisimple {{Hopf}} Actions on Regular Graded Algebras, {{I}}},
  author = {Chan, Kenneth and Kirkman, Ellen Elizabeth and Walton, Chelsea M. and Zhang, James Jian},
  date = {2018-08-15},
  journaltitle = {Journal of Algebra},
  shortjournal = {J. Algebra},
  volume = {508},
  pages = {512--538},
  doi = {10.1016/j.jalgebra.2018.05.008},
  url = {https://www.sciencedirect.com/science/article/pii/S0021869318303016},
  urldate = {2023-05-26},
  langid = {english},
  mrnumber = {3810305}
}

@article{chanMcKayCorrespondenceSemisimple2019,
  title = {{{McKay}} Correspondence for Semisimple {{Hopf}} Actions on Regular Graded Algebras. {{II}}},
  author = {Chan, Kenneth and Kirkman, Ellen Elizabeth and Walton, Chelsea M. and Zhang, James Jian},
  date = {2019-03-11},
  journaltitle = {Journal of Noncommutative Geometry},
  shortjournal = {J. Noncommutative Geom.},
  volume = {13},
  number = {1},
  pages = {87--114},
  doi = {10.4171/jncg/305},
  url = {https://ems.press/journals/jncg/articles/16076},
  urldate = {2025-07-17},
  langid = {english},
  mrnumber = {3941474}
}

@article{dlabEigenvaluesCoxeterTransformations1981,
  title = {Eigenvalues of {{Coxeter Transformations}} and the {{Gelfand-Kirillov Dimension}} of the {{Preprojective Algebras}}},
  author = {Dlab, Vlastimil B. and Ringel, Claus Michael},
  date = {1981},
  journaltitle = {Proceedings of the American Mathematical Society},
  shortjournal = {Proc. Am. Math. Soc.},
  volume = {83},
  number = {2},
  pages = {228--232},
  publisher = {American Mathematical Society},
  doi = {10.2307/2043500},
  url = {https://www.jstor.org/stable/2043500},
  urldate = {2026-01-08},
  mrnumber = {0624903}
}

@article{heCliffordDeformationsKoszul2024,
  title = {Clifford {{Deformations}} of {{Koszul Frobenius Algebras}} and {{Noncommutative Quadrics}}},
  author = {He, Ji-Wei and Ye, Yu},
  date = {2024-03},
  journaltitle = {Algebra Colloquium},
  shortjournal = {Algebra Colloq.},
  volume = {31},
  number = {1},
  pages = {63--82},
  publisher = {World Scientific Publishing Co.},
  doi = {10.1142/S1005386724000087},
  url = {https://www.worldscientific.com/doi/10.1142/S1005386724000087},
  urldate = {2025-10-16},
  mrnumber = {4717568}
}

@article{hePreresolutionsNoncommutativeIsolated2022,
  title = {Preresolutions of Noncommutative Isolated Singularities},
  author = {He, Ji-Wei and Ye, Yu},
  date = {2022-04-06},
  journaltitle = {Pacific Journal of Mathematics},
  shortjournal = {Pac. J. Math.},
  volume = {316},
  number = {2},
  pages = {367--394},
  publisher = {Mathematical Sciences Publishers},
  doi = {10.2140/pjm.2022.316.367},
  url = {https://msp.org/pjm/2022/316-2/p06.xhtml},
  urldate = {2024-09-06},
  langid = {english},
  mrnumber = {4404004}
}

@article{herschend$n$representationInfiniteAlgebras2014,
  author = {Herschend, Martin and Iyama, Osamu and Oppermann, Steffen},
  date = {2014-02-15},
  journaltitle = {Advances in Mathematics},
  shortjournal = {Adv. Math.},
  volume = {252},
  pages = {292--342},
  doi = {10.1016/j.aim.2013.09.023},
  url = {https://www.sciencedirect.com/science/article/pii/S0001870813003642},
  urldate = {2024-04-26},
  langid = {english},
  mrnumber = {3144232},
  title = {{$n$}-Representation Infinite Algebras}
}

@article{higashitaniCombinatorialStudyStable2022,
  title = {Combinatorial Study of Stable Categories of Graded {{Cohen}}--{{Macaulay}} Modules over Skew Quadric Hypersurfaces},
  author = {Higashitani, Akihiro and Ueyama, Kenta},
  date = {2022-01-01},
  journaltitle = {Collectanea Mathematica},
  shortjournal = {Collect. Math.},
  volume = {73},
  number = {1},
  pages = {43--54},
  doi = {10.1007/s13348-020-00306-1},
  url = {https://doi.org/10.1007/s13348-020-00306-1},
  urldate = {2026-10-15},
  langid = {english},
  mrnumber = {4358245}
}

@online{huClassificationNoncommutativeCentral2026,
  title = {Classification of Noncommutative Central Conics},
  author = {Hu, Hai-Gang and Mori, Izuru and Wu, Wen-Chao},
  date = {2026-02-03},
  eprint = {2602.03236v1},
  eprinttype = {arXiv},
  eprintclass = {math},
  doi = {10.48550/arXiv.2602.03236v1},
  url = {http://arxiv.org/abs/2602.03236v1},
  urldate = {2026-02-10},
  pubstate = {prepublished}
}

@article{huCliffordQuadraticComplete2024,
  title = {Clifford Quadratic Complete Intersections},
  author = {Hu, Hai-Gang and Mori, Izuru},
  date = {2024-09-04},
  journaltitle = {Mathematische Zeitschrift},
  shortjournal = {Math. Z.},
  volume = {308},
  number = {2},
  pages = {30},
  doi = {10.1007/s00209-024-03575-9},
  url = {https://doi.org/10.1007/s00209-024-03575-9},
  urldate = {2026-11-15},
  langid = {english},
  mrnumber = {4794245}
}

@article{huNoncommutativeAffinePencils2026,
  title = {Noncommutative Affine Pencils of Conics},
  author = {Hu, Hai-Gang and Mori, Izuru and Takeda, Koki and Wu, Wen-Chao},
  date = {2026-04-20},
  journaltitle = {Journal of the Australian Mathematical Society},
  shortjournal = {J. Aust. Math. Soc.},
  pages = {1--22},
  doi = {10.1017/S1446788726101530},
  url = {https://www.cambridge.org/core/journals/journal-of-the-australian-mathematical-society/article/noncommutative-affine-pencils-of-conics/65D19910CAC44D2628CB62EBA98FC6E3},
  urldate = {2026-10-21},
  langid = {english}
}

@article{huNoncommutativeConicsCalabi2023,
  title = {Noncommutative Conics in {{Calabi}}--{{Yau}} Quantum Projective Planes},
  author = {Hu, Hai-Gang and Matsuno, Masaki and Mori, Izuru},
  date = {2023-04-15},
  journaltitle = {Journal of Algebra},
  shortjournal = {J. Algebra},
  volume = {620},
  pages = {194--224},
  doi = {10.1016/j.jalgebra.2022.12.022},
  url = {https://www.sciencedirect.com/science/article/pii/S0021869322006044},
  urldate = {2025-09-19},
  mrnumber = {4531545}
}

@article{itabaASREGULARITYGEOMETRICALGEBRAS2022,
  title = {{{AS-Regularity}} of Geometric Algebras of Plane Cubic Curves},
  author = {Itaba, Ayako and Matsuno, Masaki},
  date = {2022-04},
  journaltitle = {Journal of the Australian Mathematical Society},
  shortjournal = {J. Aust. Math. Soc.},
  volume = {112},
  number = {2},
  pages = {193--217},
  doi = {10.1017/S1446788721000070},
  url = {https://www.cambridge.org/core/journals/journal-of-the-australian-mathematical-society/article/asregularity-of-geometric-algebras-of-plane-cubic-curves/C96A9B1A1660BB93D7991441D68A7DF2},
  urldate = {2026-10-20},
  langid = {english},
  mrnumber = {4389808}
}

@article{itabaDefiningRelations3dimensional2021,
  title = {Defining Relations of 3-Dimensional Quadratic {{AS}}-Regular Algebras},
  author = {Itaba, Ayako and Matsuno, Masaki},
  date = {2021},
  journaltitle = {Mathematical Journal of Okayama University},
  shortjournal = {Math. J. Okayama Univ.},
  volume = {63},
  pages = {61--86},
  mrnumber = {4174859}
}

@article{iyamaAuslanderCorrespondence2007,
  title = {Auslander Correspondence},
  author = {Iyama, Osamu},
  date = {2007-03-20},
  journaltitle = {Advances in Mathematics},
  shortjournal = {Adv. Math.},
  volume = {210},
  number = {1},
  pages = {51--82},
  doi = {10.1016/j.aim.2006.06.003},
  url = {https://www.sciencedirect.com/science/article/pii/S0001870806001733},
  urldate = {2023-09-15},
  mrnumber = {2298820}
}

@online{iyamaCohenMacaulayRepresentationsArtinSchelter2025,
  title = {Cohen--{{Macaulay}} Representations of {{Artin}}--{{Schelter Gorenstein}} Algebras of Dimension One},
  author = {Iyama, Osamu and Kimura, Yuta and Ueyama, Kenta},
  date = {2026-01-19},
  eprint = {2404.05925v4},
  eprinttype = {arXiv},
  doi = {10.48550/arXiv.2404.05925},
  url = {http://arxiv.org/abs/2404.05925v4},
  urldate = {2024-12-01},
  pubstate = {prepublished}
}

@article{iyamaFominZelevinskyMutationTilting2008,
  title = {Fomin--{{Zelevinsky}} Mutation and Tilting Modules over {{Calabi}}--{{Yau}} Algebras},
  author = {Iyama, Osamu and Reiten, Idun},
  date = {2008},
  journaltitle = {American Journal of Mathematics},
  shortjournal = {Am. J. Math.},
  volume = {130},
  number = {4},
  pages = {1087--1149},
  publisher = {Johns Hopkins University Press},
  url = {https://muse.jhu.edu/pub/1/article/244171},
  urldate = {2024-02-20},
  mrnumber = {2427009}
}

@article{iyamaMaximalModificationsAuslander2014,
  title = {Maximal Modifications and {{Auslander}}--{{Reiten}} Duality for Non-Isolated Singularities},
  author = {Iyama, Osamu and Wemyss, Michael},
  date = {2014-09-01},
  journaltitle = {Inventiones mathematicae},
  shortjournal = {Invent. Math.},
  volume = {197},
  number = {3},
  pages = {521--586},
  doi = {10.1007/s00222-013-0491-y},
  url = {https://doi.org/10.1007/s00222-013-0491-y},
  urldate = {2023-11-10},
  langid = {english},
  mrnumber = {3251829}
}

@inproceedings{iyamaTiltingCohenMacaulayRepresentations2018,
  title = {Tilting {{Cohen-Macaulay}} Representations},
  booktitle = {Proceedings of the {{International Congress}} of {{Mathematicians}}},
  author = {Iyama, Osamu},
  date = {2018},
  volume = {II},
  pages = {125--162},
  publisher = {World Scientific},
  location = {Hackensack, New Jersey},
  eventdate = {2018-08-01/2018-08-09},
  eventtitle = {International {{Congress}} of {{Mathematicians}}},
  mrnumber = {3966760},
  venue = {Rio de Janeiro, Brazil}
}

@article{iyaTiltingClusterTilting2013,
  title = {Tilting and Cluster Tilting for Quotient Singularities},
  author = {Iyama, Osamu and Takahashi, Ryo},
  date = {2013-07-01},
  journaltitle = {Mathematische Annalen},
  shortjournal = {Math. Ann.},
  volume = {356},
  number = {3},
  pages = {1065--1105},
  doi = {10.1007/s00208-012-0842-9},
  url = {https://doi.org/10.1007/s00208-012-0842-9},
  urldate = {2023-07-25},
  langid = {english},
  mrnumber = {3063907}
}

@article{kapranovDerivedCategoriesCoherent1988,
  title = {On the Derived Categories of Coherent Sheaves on Some Homogeneous Spaces},
  author = {Kapranov, Mikhail M.},
  date = {1988-10-01},
  journaltitle = {Inventiones mathematicae},
  shortjournal = {Invent. Math.},
  volume = {92},
  number = {3},
  pages = {479--508},
  doi = {10.1007/BF01393744},
  url = {https://doi.org/10.1007/BF01393744},
  urldate = {2025-12-24},
  langid = {english},
  mrnumber = {0939472}
}

@article{kapranovDerivedCategoryCoherent1986,
  title = {The Derived Category of Coherent Sheaves on the Square},
  author = {Kapranov, Mikhail M.},
  date = {1986-04-01},
  journaltitle = {Functional Analysis and Its Applications},
  shortjournal = {Funct. Anal. Its Appl.},
  volume = {20},
  number = {2},
  pages = {141--142},
  doi = {10.1007/BF01077274},
  url = {https://doi.org/10.1007/BF01077274},
  urldate = {2025-12-24},
  langid = {english},
  mrnumber = {0847146}
}

@article{kellerDeformedCalabiYau2011,
  title = {Deformed {{Calabi}}--{{Yau}} Completions},
  author = {Keller, Bernhard},
  namea = {family=Bergh, given=Michel, prefix=van den, useprefix=true},
  nameatype = {collaborator},
  date = {2011-05-01},
  journaltitle = {Journal f\"ur die reine und angewandte Mathematik},
  shortjournal = {J. Reine Angew. Math.},
  volume = {2011},
  number = {654},
  pages = {125--180},
  publisher = {De Gruyter},
  doi = {10.1515/crelle.2011.031},
  url = {https://www.degruyterbrill.com/document/doi/10.1515/crelle.2011.031/html},
  urldate = {2026-03-30},
  langid = {english},
  mrnumber = {2795754},
  note = {With an appendix by Michel van den Bergh}
}

@book{lamFirstCourseNoncommutative2001,
  title = {A {{First Course}} in {{Noncommutative Rings}}},
  author = {Lam, Tsit Yuen},
  date = {2001},
  series = {Graduate {{Texts}} in {{Mathematics}}},
  edition = {2},
  number = {131},
  publisher = {Springer-Verlag},
  location = {New York, New York},
  doi = {10.1007/978-1-4419-8616-0},
  url = {http://link.springer.com/10.1007/978-1-4419-8616-0},
  urldate = {2023-09-19},
  mrnumber = {1838439},
  pagetotal = {xx+385 pp.},
  shortseries = {Grad. Texts in Math.}
}

@inproceedings{lenzingHomologicalTransferFinitely1983,
  title = {Homological {{Transfer}} from {{Finitely Presented}} to {{Infinite Modules}}},
  booktitle = {Abelian {{Group Theory}}},
  author = {Lenzing, Helmut},
  editor = {G\"obel, R\"udiger and Lady, Lee and Mader, Adolf},
  date = {1983},
  pages = {734--761},
  publisher = {Springer},
  location = {Berlin, Germany},
  doi = {10.1007/978-3-662-21560-9_53},
  url = {https://doi.org/10.1007/978-3-662-21560-9_53},
  urldate = {2026-03-30},
  eventdate = {1982-12-18/1983-01-04},
  eventtitle = {The Conference Held at the {{University}} of {{Hawaii}}},
  langid = {english},
  mrnumber = {0722664},
  venue = {Honolulu, Hawaii}
}

@article{leuschkeEndomorphismRingsFinite2007,
  title = {Endomorphism {{Rings}} of {{Finite Global Dimension}}},
  author = {Leuschke, Graham J.},
  date = {2007-04},
  journaltitle = {Canadian Journal of Mathematics},
  shortjournal = {Can. J. Math.},
  volume = {59},
  number = {2},
  pages = {332--342},
  doi = {10.4153/CJM-2007-014-1},
  url = {https://www.cambridge.org/core/journals/canadian-journal-of-mathematics/article/endomorphism-rings-of-finite-global-dimension/39CDE82CCE2E87D9B082F9D7C9D7FAB3},
  urldate = {2025-07-17},
  langid = {english},
  mrnumber = {2310620}
}

@article{levasseurPropertiesNoncommutativeRegular1992,
  title = {Some Properties of Non-Commutative Regular Graded Rings},
  author = {Levasseur, Thierry},
  date = {1992-09},
  journaltitle = {Glasgow Mathematical Journal},
  shortjournal = {Glasg. Math. J.},
  volume = {34},
  number = {3},
  pages = {277--300},
  publisher = {Cambridge University Press},
  doi = {10.1017/S0017089500008843},
  url = {https://www.cambridge.org/core/journals/glasgow-mathematical-journal/article/some-properties-of-noncommutative-regular-graded-rings/244AE0A8E9C8AB6782515B504F1AC5C0#},
  urldate = {2023-03-06},
  langid = {english},
  mrnumber = {1181768}
}

@online{liCommonlyGradedAlgebras2025,
  title = {Commonly {{Graded Algebras}} and {{Their Homological Properties}}},
  author = {Li, Hao-Nan and Wu, Quan-Shui},
  date = {2025-08-08},
  eprint = {2508.06068v1},
  eprinttype = {arXiv},
  eprintclass = {math},
  doi = {10.48550/arXiv.2508.06068v1},
  url = {http://arxiv.org/abs/2508.06068v1},
  urldate = {2025-10-14},
  pubstate = {prepublished}
}

@article{liGeneralizedKoszulAlgebra2023,
  title = {Generalized {{Koszul}} Algebra and {{Koszul}} Duality},
  author = {Li, Hao-Nan and Wu, Quan-Shui},
  date = {2023-04-01},
  journaltitle = {Journal of Algebra},
  shortjournal = {J. Algebra},
  volume = {619},
  pages = {643--685},
  doi = {10.1016/j.jalgebra.2022.12.023},
  url = {https://www.sciencedirect.com/science/article/pii/S0021869322006020},
  urldate = {2023-11-30},
  mrnumber = {4530820}
}

@article{liNoncommutativeResolutionsASGorenstein2026,
  title = {Noncommutative Resolutions of {{AS-Gorenstein}} Isolated Singularities},
  author = {Li, Hao-Nan and Shen, Men-Da and Wu, Quan-Shui},
  date = {2026-06-15},
  journaltitle = {Journal of Noncommutative Geometry},
  shortjournal = {J. Noncommutative Geom.},
  doi = {10.4171/jncg/670},
  url = {https://ems.press/journals/jncg/articles/14299817},
  urldate = {2026-09-30},
  langid = {english},
  publication-number = {670},
  note = {online first}
}

@online{liNoncommutativeResolutionsNoncommutative2026,
  title = {Noncommutative Resolutions of Noncommutative Isolated Singularities},
  author = {Li, Hao-Nan and Wu, Quan-Shui},
  date = {2026-09-18},
  eprint = {2609.21409},
  eprinttype = {arXiv},
  eprintclass = {math.RA},
  doi = {10.48550/arXiv.2609.21409},
  url = {http://arxiv.org/abs/2609.21409},
  urldate = {2026-09-24},
  pubstate = {prepublished}
}

@article{liuSkewKnorrersPeriodicity2026,
  title = {Skew {{Kn\"orrer}}'s Periodicity Theorem},
  author = {Liu, Yang and Shen, Yuan and Wang, Xin},
  date = {2026-01-20},
  journaltitle = {Journal of Noncommutative Geometry},
  shortjournal = {J. Noncommutative Geom.},
  volume = {20},
  number = {4},
  pages = {1465--1501},
  doi = {10.4171/jncg/652},
  url = {https://ems.press/journals/jncg/articles/14299442},
  urldate = {2026-10-30},
  langid = {english},
  mrnumber = {5130084}
}

@inproceedings{martinez-villaApplicationsKoszulAlgebras1996,
  title = {Applications of {{Koszul}} Algebras: The Preprojective Algebra},
  shorttitle = {Applications of {{Koszul}} Algebras},
  booktitle = {Representation Theory of Algebras ({{Cocoyoc}}, 1994)},
  author = {Mart\'inez-Villa, Roberto},
  date = {1996},
  series = {{{CMS Conf}}. {{Proc}}.},
  volume = {18},
  pages = {487--504},
  publisher = {Amer. Math. Soc., Providence, RI},
  mrnumber = {1388069}
}

@article{martinez-villaSkewGroupAlgebras2001,
  title = {Skew Group Algebras and Their {{Yoneda}} Algebras},
  author = {Mart\'inez-Villa, Roberto},
  date = {2001},
  journaltitle = {Mathematical Journal of Okayama University},
  shortjournal = {Math. J. Okayama Univ.},
  volume = {43},
  pages = {1--16},
  mrnumber = {1913868}
}

@article{minamotoAmplenessTwosidedTilting2012,
  title = {Ampleness of Two-Sided Tilting Complexes},
  author = {Minamoto, Hiroyuki},
  date = {2012-01-01},
  journaltitle = {International Mathematics Research Notices},
  shortjournal = {Int. Math. Res. Not.},
  volume = {2012},
  number = {1},
  pages = {67--101},
  doi = {10.1093/imrn/rnr001},
  url = {https://doi.org/10.1093/imrn/rnr001},
  urldate = {2024-03-28},
  langid = {english},
  mrnumber = {2874928}
}

@article{minamotoStructureASGorensteinAlgebras2011,
  title = {The Structure of {{AS-Gorenstein}} Algebras},
  author = {Minamoto, Hiroyuki and Mori, Izuru},
  date = {2011-03-20},
  journaltitle = {Advances in Mathematics},
  shortjournal = {Adv. Math.},
  volume = {226},
  number = {5},
  pages = {4061--4095},
  doi = {10.1016/j.aim.2010.11.004},
  url = {https://www.sciencedirect.com/science/article/pii/S0001870810004226},
  urldate = {2023-03-22},
  langid = {english},
  mrnumber = {2770441}
}

@article{mitchellRingsSeveralObjects1972,
  title = {Rings with Several Objects},
  author = {Mitchell, Barry},
  date = {1972-02-01},
  journaltitle = {Advances in Mathematics},
  shortjournal = {Adv. Math.},
  volume = {8},
  number = {1},
  pages = {1--161},
  doi = {10.1016/0001-8708(72)90002-3},
  url = {https://www.sciencedirect.com/science/article/pii/0001870872900023},
  urldate = {2026-04-01},
  mrnumber = {0294454}
}

@article{moriAmpleGroupAction2016,
  title = {Ample Group Action on {{AS}}-Regular Algebras and Noncommutative Graded Isolated Singularities},
  author = {Mori, Izuru and Ueyama, Kenta},
  date = {2016-10},
  journaltitle = {Transactions of the American Mathematical Society},
  shortjournal = {Trans. Am. Math. Soc.},
  volume = {368},
  number = {10},
  pages = {7359--7383},
  doi = {10.1090/tran/6580},
  url = {https://www.ams.org/tran/2016-368-10/S0002-9947-2015-06580-5/},
  urldate = {2023-05-24},
  langid = {english},
  mrnumber = {3471094}
}

@article{moriCategoricalCharacterizationQuantum2021,
  title = {A Categorical Characterization of Quantum Projective Spaces},
  author = {Mori, Izuru and Ueyama, Kenta},
  date = {2021-09-01},
  journaltitle = {Journal of Noncommutative Geometry},
  shortjournal = {J. Noncommutative Geom.},
  volume = {15},
  number = {2},
  pages = {489--529},
  doi = {10.4171/jncg/403},
  url = {https://ems.press/journals/jncg/articles/2504083},
  urldate = {2024-03-17},
  langid = {english}
}

@article{moriHomologicalPropertiesBalanced2003,
  title = {Homological Properties of Balanced {{Cohen}}--{{Macaulay}} Algebras},
  author = {Mori, Izuru},
  date = {2003},
  journaltitle = {Transactions of the American Mathematical Society},
  shortjournal = {Trans. Am. Math. Soc.},
  volume = {355},
  number = {3},
  pages = {1025--1042},
  doi = {10.1090/S0002-9947-02-03166-5},
  url = {https://www.ams.org/tran/2003-355-03/S0002-9947-02-03166-5/},
  urldate = {2023-06-09},
  langid = {english},
  mrnumber = {1938744}
}

@online{moriIrreducibleNoncommutativeQuadrics2026,
  title = {Irreducible Noncommutative Quadrics},
  author = {Mori, Izuru and Ueyama, Kenta and Wu, Wen-Chao},
  date = {2026-06-23},
  eprint = {2606.24606v1},
  eprinttype = {arXiv},
  eprintclass = {math.RA},
  doi = {10.48550/arXiv.2606.24606},
  url = {http://arxiv.org/abs/2606.24606v1},
  urldate = {2026-10-19},
  pubstate = {prepublished}
}

@article{moriMcKaytypeCorrespondenceASregular2013,
  title = {{{McKay}}-Type Correspondence for {{AS}}-Regular Algebras},
  author = {Mori, Izuru},
  date = {2013},
  journaltitle = {Journal of the London Mathematical Society},
  shortjournal = {J. Lond. Math. Soc.},
  volume = {88},
  number = {1},
  pages = {97--117},
  doi = {10.1112/jlms/jdt005},
  url = {https://onlinelibrary.wiley.com/doi/abs/10.1112/jlms/jdt005},
  urldate = {2026-03-12},
  langid = {english},
  mrnumber = {3092260}
}

@article{moriNoncommutativeKnorrersPeriodicity2022,
  title = {Noncommutative {{Kn\"orrer}}'s Periodicity Theorem and Noncommutative Quadric Hypersurfaces},
  author = {Mori, Izuru and Ueyama, Kenta},
  date = {2022-04-27},
  journaltitle = {Algebra \& Number Theory},
  shortjournal = {Algebra Number Theory},
  volume = {16},
  number = {2},
  pages = {467--504},
  publisher = {Mathematical Sciences Publishers},
  doi = {10.2140/ant.2022.16.467},
  url = {https://msp.org/ant/2022/16-2/p09.xhtml},
  urldate = {2025-02-17},
  mrnumber = {4412580}
}

@article{moriNoncommutativeMatrixFactorizations2021,
  title = {Noncommutative Matrix Factorizations with an Application to Skew Exterior Algebras},
  author = {Mori, Izuru and Ueyama, Kenta},
  date = {2021-11-15},
  journaltitle = {Journal of Algebra},
  shortjournal = {J. Algebra},
  volume = {586},
  pages = {1053--1087},
  doi = {10.1016/j.jalgebra.2021.07.012},
  url = {https://www.sciencedirect.com/science/article/pii/S0021869321003690},
  urldate = {2025-12-26},
  mrnumber = {4296228}
}

@article{morStableCategoriesGraded2016,
  title = {Stable Categories of Graded Maximal {{Cohen}}--{{Macaulay}} Modules over Noncommutative Quotient Singularities},
  author = {Mori, Izuru and Ueyama, Kenta},
  date = {2016-07-16},
  journaltitle = {Advances in Mathematics},
  shortjournal = {Adv. Math.},
  volume = {297},
  pages = {54--92},
  doi = {10.1016/j.aim.2016.04.009},
  url = {https://www.sciencedirect.com/science/article/pii/S0001870815302103},
  urldate = {2023-06-02},
  langid = {english},
  mrnumber = {3498794}
}

@incollection{orlovDerivedCategoriesCoherent2009,
  title = {Derived {{Categories}} of {{Coherent Sheaves}} and {{Triangulated Categories}} of {{Singularities}}},
  booktitle = {Algebra, {{Arithmetic}}, and {{Geometry}}},
  author = {Orlov, Dmitri},
  editor = {Tschinkel, Yuri and Zarhin, Yuri},
  date = {2009},
  series = {Progress in {{Mathematics}}},
  volume = {2},
  number = {269},
  pages = {503--531},
  publisher = {Birkh\"auser},
  location = {Boston, Massachusetts},
  doi = {10.1007/978-0-8176-4747-6_16},
  url = {https://doi.org/10.1007/978-0-8176-4747-6_16},
  urldate = {2023-11-16},
  langid = {english},
  mrnumber = {2641200},
  shortseries = {Progr. Math.},
  volumes = {2}
}

@article{qinNoncommutativeQuasiresolutions2019,
  title = {Noncommutative Quasi-Resolutions},
  author = {Qin, Xiao-Shan and Wang, Yan-Hua and Zhang, James Jian},
  date = {2019-10-15},
  journaltitle = {Journal of Algebra},
  shortjournal = {J. Algebra},
  volume = {536},
  pages = {102--148},
  doi = {10.1016/j.jalgebra.2019.07.015},
  url = {https://www.sciencedirect.com/science/article/pii/S0021869319303904},
  urldate = {2024-09-09},
  mrnumber = {3989612}
}

@article{reitenNoetherianHereditaryAbelian2002,
  title = {Noetherian Hereditary Abelian Categories Satisfying {{Serre}} Duality},
  author = {Reiten, Idun and family=Bergh, given=Michel, prefix=van den, useprefix=true},
  date = {2002-04-01},
  journaltitle = {Journal of the American Mathematical Society},
  shortjournal = {J. Am. Math. Soc.},
  volume = {15},
  number = {2},
  pages = {295--366},
  doi = {10.1090/S0894-0347-02-00387-9},
  url = {https://www.ams.org/jams/2002-15-02/S0894-0347-02-00387-9/},
  urldate = {2023-04-28},
  langid = {english},
  mrnumber = {1887637},
  shortauthor = {RV}
}

@article{reyesGradedTwistedCalabi2022,
  title = {Graded {{Twisted Calabi}}--{{Yau Algebras Are Generalized Artin}}--{{Schelter Regular}}},
  author = {Reyes, Manuel Lionel and Rogalski, Daniel},
  date = {2022-03},
  journaltitle = {Nagoya Mathematical Journal},
  shortjournal = {Nagoya Math. J.},
  volume = {245},
  pages = {100--153},
  publisher = {Cambridge University Press},
  doi = {10.1017/nmj.2020.32},
  url = {https://www.cambridge.org/core/journals/nagoya-mathematical-journal/article/graded-twisted-calabiyau-algebras-are-generalized-artinschelter-regular/751A2B1B217F110A5EC754710940ED3A},
  urldate = {2022-10-04},
  langid = {english},
  mrnumber = {4413364}
}

@article{reyesGrowthGradedTwisted2019,
  title = {Growth of Graded Twisted {{Calabi-Yau}} Algebras},
  author = {Reyes, Manuel L. and Rogalski, Daniel},
  date = {2019-12-01},
  journaltitle = {Journal of Algebra},
  shortjournal = {J. Algebra},
  volume = {539},
  pages = {201--259},
  doi = {10.1016/j.jalgebra.2019.07.029},
  url = {https://www.sciencedirect.com/science/article/pii/S0021869319304338},
  urldate = {2025-03-03},
  mrnumber = {3995777}
}

@article{reyesSkewCalabiYau2014,
  title = {Skew {{Calabi}}--{{Yau}} Algebras and Homological Identities},
  author = {Reyes, Manuel and Rogalski, Daniel and Zhang, James J.},
  date = {2014-10-20},
  journaltitle = {Advances in Mathematics},
  shortjournal = {Adv. Math.},
  volume = {264},
  pages = {308--354},
  doi = {10.1016/j.aim.2014.07.010},
  url = {https://www.sciencedirect.com/science/article/pii/S0001870814002424},
  urldate = {2023-03-05},
  langid = {english},
  mrnumber = {3250287}
}

@article{sheltonKoszulAlgebrasNew2001,
  title = {On {{Koszul Algebras}} and a {{New Construction}} of {{Artin}}--{{Schelter Regular Algebras}}},
  author = {Shelton, Brad and Tingey, Craig},
  date = {2001-07-15},
  journaltitle = {Journal of Algebra},
  shortjournal = {J. Algebra},
  volume = {241},
  number = {2},
  pages = {789--798},
  doi = {10.1006/jabr.2001.8781},
  url = {https://www.sciencedirect.com/science/article/pii/S0021869301987812},
  urldate = {2025-10-16},
  mrnumber = {1843325}
}

@inproceedings{smithFinitedimensionalAlgebrasRelated1996,
  title = {Some Finite-Dimensional Algebras Related to Elliptic Curves},
  booktitle = {Representation Theory of Algebras and Related Topics},
  author = {Smith, S. Paul},
  editor = {Bautista, Raymundo and Mart\'inez-Villa, Roberto and family=Pe\~na, given=Jos\'e Antonio, prefix=de la, useprefix=true},
  date = {1996},
  series = {Canadian {{Mathematical Society Conference Proceedings}}},
  number = {19},
  pages = {315--348},
  publisher = {Amer. Math. Soc.},
  location = {Providence, Rhode Island},
  eventdate = {1994-08-16/1994-08-20},
  eventtitle = {Workshop on {{Representation Theory}} of {{Algebras}} and {{Related Topics}}},
  mrnumber = {1388568},
  shortseries = {CMS Conf. Proc.},
  venue = {Mexico City, Mexico}
}

@article{smithNoncommutativeQuadricSurfaces2013,
  title = {Noncommutative Quadric Surfaces},
  author = {Smith, S. Paul and family=Bergh, given=Michel, prefix=van den, useprefix=true},
  date = {2013-09-24},
  journaltitle = {Journal of Noncommutative Geometry},
  shortjournal = {J. Noncommutative Geom.},
  volume = {7},
  number = {3},
  pages = {817--856},
  doi = {10.4171/jncg/136},
  url = {https://ems.press/journals/jncg/articles/4781},
  urldate = {2024-04-23},
  langid = {english},
  mrnumber = {3108697},
  shortauthor = {SV}
}

@article{stephensonGrowthGradedNoetherian1997,
  title = {Growth of Graded Noetherian Rings},
  author = {Stephenson, Darin R. and Zhang, James Jian},
  date = {1997},
  journaltitle = {Proceedings of the American Mathematical Society},
  shortjournal = {Proc. Am. Math. Soc.},
  volume = {125},
  number = {6},
  pages = {1593--1605},
  doi = {10.1090/S0002-9939-97-03752-0},
  url = {https://www.ams.org/proc/1997-125-06/S0002-9939-97-03752-0/},
  urldate = {2025-03-05},
  langid = {english},
  mrnumber = {1371143}
}

@article{ueyamaClusterTiltingModules2017,
  title = {Cluster Tilting Modules and Noncommutative Projective Schemes},
  author = {Ueyama, Kenta},
  date = {2017-06-19},
  journaltitle = {Pacific Journal of Mathematics},
  shortjournal = {Pac. J. Math.},
  volume = {289},
  number = {2},
  pages = {449--468},
  publisher = {Mathematical Sciences Publishers},
  doi = {10.2140/pjm.2017.289.449},
  url = {https://msp.org/pjm/2017/289-2/p08.xhtml},
  urldate = {2023-07-26},
  langid = {english},
  mrnumber = {3667180}
}

@article{ueyamaDerivedCategoriesSkew2023,
  title = {Derived Categories of Skew Quadric Hypersurfaces},
  author = {Ueyama, Kenta},
  date = {2023-03-01},
  journaltitle = {Israel Journal of Mathematics},
  shortjournal = {Isr. J. Math.},
  volume = {253},
  number = {1},
  pages = {205--247},
  doi = {10.1007/s11856-022-2360-0},
  url = {https://doi.org/10.1007/s11856-022-2360-0},
  urldate = {2024-06-07},
  langid = {english},
  mrnumber = {4575408}
}

@article{ueyamaGradedMaximalCohen2013,
  title = {Graded Maximal {{Cohen}}--{{Macaulay}} Modules over Noncommutative Graded {{Gorenstein}} Isolated Singularities},
  author = {Ueyama, Kenta},
  date = {2013-06-01},
  journaltitle = {Journal of Algebra},
  shortjournal = {J. Algebra},
  volume = {383},
  pages = {85--103},
  doi = {10.1016/j.jalgebra.2013.02.022},
  url = {https://www.sciencedirect.com/science/article/pii/S0021869313001105},
  urldate = {2023-04-28},
  langid = {english},
  mrnumber = {3037969}
}

@article{ueyamaKnorrerPeriodicityQuadric2019,
  title = {On {{Kn\"orrer Periodicity}} for {{Quadric Hypersurfaces}} in {{Skew Projective Spaces}}},
  author = {Ueyama, Kenta},
  date = {2019-12},
  journaltitle = {Canadian Mathematical Bulletin},
  shortjournal = {Can. Math. Bull.},
  volume = {62},
  number = {4},
  pages = {896--911},
  doi = {10.4153/S0008439518000607},
  url = {https://www.cambridge.org/core/journals/canadian-mathematical-bulletin/article/on-knorrer-periodicity-for-quadric-hypersurfaces-in-skew-projective-spaces/90BDF8E9B0B6A8AE71247A49F9A97334},
  urldate = {2026-05-05},
  langid = {english},
  mrnumber = {4028496}
}

@article{ueyamaNoncommutativeGradedAlgebras2015,
  title = {Noncommutative Graded Algebras of Finite {{Cohen}}--{{Macaulay}} Representation Type},
  author = {Ueyama, Kenta},
  date = {2015-09},
  journaltitle = {Proceedings of the American Mathematical Society},
  shortjournal = {Proc. Am. Math. Soc.},
  volume = {143},
  number = {9},
  pages = {3703--3715},
  doi = {10.1090/proc/12527},
  url = {https://www.ams.org/proc/2015-143-09/S0002-9939-2015-12527-0/},
  urldate = {2023-05-31},
  langid = {english},
  mrnumber = {3359563}
}

@inproceedings{vandenberghNoncommutativeCrepantResolutions2004,
  title = {Non-Commutative {{Crepant Resolutions}}},
  booktitle = {The {{Legacy}} of {{Niels Henrik Abel}}},
  author = {family=Bergh, given=Michel, prefix=van den, useprefix=true},
  editor = {Laudal, Olav Arnfinn and Piene, Ragni},
  date = {2004},
  pages = {749--770},
  publisher = {Springer},
  location = {Berlin, Heidelberg},
  doi = {10.1007/978-3-642-18908-1_26},
  url = {https://doi.org/10.1007/978-3-642-18908-1_26},
  urldate = {2024-01-02},
  eventdate = {2002},
  eventtitle = {Abel {{Bicentennial Conference}}},
  langid = {english},
  mrnumber = {2077594},
  shortauthor = {VdB},
  venue = {Oslo}
}

@article{yekDualizingComplexesNoncommutative1992,
  title = {Dualizing {{Complexes}} over {{Noncommutative Graded Algebras}}},
  author = {Yekutieli, Amnon},
  date = {1992-12-01},
  journaltitle = {Journal of Algebra},
  shortjournal = {J. Algebra},
  volume = {153},
  number = {1},
  pages = {41--84},
  doi = {10.1016/0021-8693(92)90148-F},
  url = {https://www.sciencedirect.com/science/article/pii/002186939290148F},
  urldate = {2022-10-06},
  langid = {english},
  mrnumber = {1195406}
}

@article{yekutieliSerreDualityNoncommutative1997,
  title = {Serre Duality for Noncommutative Projective Schemes},
  author = {Yekutieli, Amnon and Zhang, James Jian},
  date = {1997},
  journaltitle = {Proceedings of the American Mathematical Society},
  shortjournal = {Proc. Am. Math. Soc.},
  volume = {125},
  number = {3},
  pages = {697--707},
  doi = {10.1090/S0002-9939-97-03782-9},
  url = {https://www.ams.org/proc/1997-125-03/S0002-9939-97-03782-9/},
  urldate = {2023-04-28},
  langid = {english},
  mrnumber = {1372045}
}

@article{zhangTwistedGradedAlgebras1996,
  title = {Twisted {{Graded Algebras}} and {{Equivalences}} of {{Graded Categories}}},
  author = {Zhang, James Jian},
  date = {1996},
  journaltitle = {Proceedings of the London Mathematical Society. Third Series},
  shortjournal = {Proc. Lond. Math. Soc. Third Ser.},
  volume = {72},
  number = {2},
  pages = {281--311},
  doi = {10.1112/plms/s3-72.2.281},
  url = {https://onlinelibrary.wiley.com/doi/abs/10.1112/plms/s3-72.2.281},
  urldate = {2023-11-29},
  langid = {english},
  mrnumber = {1367080}
}
\end{document}